\documentclass{article}
\usepackage{amssymb,amsmath}
\usepackage{url}
\usepackage{hyperref}
\usepackage{srcltx}
\usepackage{xcolor}

\def\bU{{\mathbf{U}}}

\def\bEx{{\mathbf{E}}}
\def\VI{{\mathrm{VI}}}
\def\epsVI{{\epsilon_{\hbox{\tiny\rm VI}}}}

\def\three?{3}
\def\four?{4}
\def\ten?{10}
\def\bM{{\mathbf{M}}}
\def\bK{{\mathbf{K}}}
\def\bL{{\mathbf{L}}}

\newtheorem{remark}{Remark}[section]
\newtheorem{proposition}{Proposition}[section]

\def\norm2to2{{\|\cdot\|_{2,2}}}

\def\cl{{\hbox{\rm  cl}\,}}

\def\inter{\hbox{\rm  int}}
\def\K{{\mathbf{K}}}

\def\Diag{\hbox{\rm  Diag}}

\def\K{{\mathbf{K}}}

\def\Conv{\hbox{\rm  Conv}}

\def\Tr{{\mathop{\hbox{\rm  Tr}}}}
\def\cA{{\cal A}}

\def\cC{{\cal C}}

\def\cF{{\cal F}}

\def\cK{{\cal K}}

\def\cP{{\cal P}}
\def\cQ{{\cal Q}}

\def\cT{{\cal T}}
\def\cU{{\cal U}}
\def\cV{{\cal V}}

\def\cX{{\cal X}}
\def\cY{{\cal Y}}
\def\cZ{{\cal Z}}

\def\F{{\cal F}}

\def\K{{\cal K}}

\def\P{{\cal P}}
\def\Q{{\cal Q}}
\def\R{{\cal R}}

\def\T{{\cal T}}
\def\U{{\cal U}}
\def\V{{\cal V}}

\def\X{{\cal X}}
\def\Y{{\cal Y}}
\def\Z{{\cal Z}}

\def\Dom{{\mathop{\hbox{\rm  Dom}\,}}}

\def\Epi{{\hbox{\rm  Epi}}}

\def\bK{{\mathbf{K}}}
\def\bL{{\mathbf{L}}}
\def\bS{{\mathbf{S}}}

\def\e{{\hbox{\rm e}}}

\def\qed{\ \hfill$\square$\par\smallskip}

\def\mypict3{\epsfxsize=220pt\epsfysize=80pt\epsfbox}

\def\bR{{\mathbf{R}}}

\newcommand{\hide}[1]{{}}

\title{On Well-Structured Convex-Concave Saddle Point Problems and Variational Inequalities with Monotone Operators}
\author{
Anatoli Juditsky
\thanks{LJK, Universit\'e Grenoble Alpes, 700 Avenue Centrale,  38401 Domaine Universitaire de Saint-Martin-d'Hères, France,	
{\tt anatoli.juditsky@univ-grenoble-alpes.fr}}
\and Arkadi Nemirovski
\thanks{ISyE, Georgia Institute
 of Technology, Atlanta, Georgia
30332, USA, {\tt nemirovs@isye.gatech.edu}.}}
\date{}
\begin{document}
\maketitle
\begin{abstract}
For those acquainted with {\tt CVX} ({aka} disciplined convex programming) of M. Grant and S. Boyd \cite{cvx}, the motivation {of this work} is  the desire to {extend the scope of} {\tt CVX} {beyond  convex minimization---to} convex-concave saddle point problems and variational inequalities with monotone operators. To {attain this goal},
given a family $\cK$ of cones (e.g., {Lorentz, }semidefinite{, geometric, etc}), we introduce the notions of $\cK$-conic representation of a convex-concave saddle point problem and of variational inequality with monotone operator. We demonstrate that given such a representation of the problem of interest, the latter can be reduced straightforwardly to a conic problem on a cone from $\cK$ and thus can be solved by  (any) solver capable to handle conic problems on cones from $\cK$ (e.g., {\tt Mosek} or {\tt SDPT3} in the case of semidefinite cones). We also  show that $\cK$-representations of convex-concave functions and monotone vector fields admit a  fully algorithmic calculus which helps to recognize the cases when a saddle point{ problem or }variational inequality  can be converted into a conic problem on a cone from $\cK$ and to carry out such conversion.
\end{abstract}
\section{Introduction}
{Along with emergence of powerful computational tools, one of important} components of what is {sometimes referred to as} Interior-Point Revolution in Convex Optimization was the rise of (informal) notion of a ``well-structured convex problem''---{convex optimization problem} which can be reformulated as a conic problem on a cone from a ``nice'' family $\cK$, most notably, the families of finite direct products of nonnegative rays (Linear Programming), Lorentz cones (Second Order Conic Programming), and semidefinite cones (Semidefinite Programming). Problems of the latter ``magic families'' cover basically all needs arising in applications of Convex Optimization to decision making, engineering, statistics, etc. At the same time, reformulating problem of interest as LP/SOCP/SDP makes the problem amenable for numerical processing by  a ``universal'' SDP solver, so that end user should not care much about number-crunching. {Today,} there exist powerful universal SDP solvers {such as} {\tt SDPT3} \cite{SDPT3} or industrial-grade {\tt Mosek} \cite{Mosek}, with permanently improving performance and reliability. \par
The ``maiden form'' of an optimization problem arising in applications usually is not
the ``nice conic'' one, even when at the end of the day the problem can be {cast into,} say, an SDP program. {Such transformation is carried out with the help of a special} ``calculus'' (see, e.g., \cite{BoydVan,LMCOBook,LMCOLN}) which allows to recognize that the problem at hand can be reformulated as, say, an SDP. This calculus, as any other, has two components: (a) ``raw materials''---a collection of sets and functions {possessing} {explicit} representations via SDP (what exactly ``representation'' means, will be explained in Section \ref{prelim}), and (b) ``calculus rules'' expressing SDP representation of the result of a specific operation with functions/sets via SDP representations of the operands. It turns out that calculus rules {for SDP representations} cover all basic convexity-preserving operations with functions (like taking linear combinations with nonnegative coefficients and convexity-preserving compositions) and sets (like taking finite intersections and images/inverse images under affine mappings). Moreover, the calculus rules happen to be fully algorithmic and thus can be implemented on a compiler. This possibility is utilized{, for instance,} in the {\tt CVX} software of M. Grant and S. Boyd \cite{cvx} which is the second-to-none in terms of its scope and user-friendliness tool for numerical processing of convex optimization problems.
\par
So far, we were speaking about (constrained) convex minimization, and this is the {usual} scope of the calculus of well-structured  convex problems. For example, given a number of SDP-representable functions, the calculus allows to obtain SDP representation of their maximum and thus convert to SDP program, in an automated fashion, {the problem of} minimizing this maximum over a feasible set given by an SDP representation. However, such calculus does {\sl not}
allow to get an SDP representation of the maximum of {\sl infinite} family of SDP-representable functions, even when this family is perfectly well organized. As a result, convex-concave saddle point problems,
well-structured from the ``viewpoint of a human,'' cannot be ``fed to CVX'' as they are, and reducing such a problem to the conic form amenable for the existing powerful  convex optimization software requires ad hoc work on case-to-case basis. Similar difficulties arise when instead of convex-concave saddle point problems {someone is looking to solve numerically} other problems with convex structure, most notably, variational inequalities (VI) with monotone operators. Monotone VI's can be thought of as the most general problems with convex structure: on one hand, other problems of this type (constrained convex minimization, convex-concave saddle points, convex Nash equilibrium problems) under extremely mild boundedness and regularity assumptions  can be reduced to monotone VI's (for details, see e.g., \cite[Section 5.6]{LMCOLN}).
On the other hand, efficient {\sl black box oriented} algorithms of convex minimization, like First Order or Ellipsoid  methods, can be extended to monotone VI's  (cf., e.g., \cite{NOR}) without {harming} their efficiency estimates. There exists also {rich} literature on extending polynomial time Interior-Point algorithms onto linear and monotone complementarity problems, see, e.g., \cite{Kojima,Renato,Pang,Potra,RW97,Sim} and references therein. In particular, it is well known that VI with {\sl affine} monotone operator and domain  admitting conic representation can be straightforwardly reduced to a conic optimization problem. In \cite[Section 7.4]{NNBook} this fact was somehow extended
to  nonlinear monotone VI's from certain restricted class---those admitting an explicit ``convex representation'' allowing to reduce the VI of interest to an affine monotone VI on a properly modified domain; that reference develops also a rudimentary calculus of convex representations of monotone VI's. {Some results\footnote{dealing with the search for strong solutions to VI's with not necessarily monotone operator $F$
under rather restrictive assumptions that $x^TF(x)$ and $-F(x)$ are convex on polyhedral set $\X\subset\bR^n_+$}  on reformulations of VI's as convex optimization problems can be found in \cite{ABP2006}.}
That being said, to the best of our knowledge, {beyond the classes of VI's with affine monotone operators and related bilinear saddle point problems} there is no {much} understanding of what is a ``well structured'' convex-concave saddle point problem or {``genuine''---not potential---}monotone VI, and of how such a problem should be represented in order for the representation to admit a meaningful algorithmic calculus, on one hand, and to allow for straightforward conversion of the problem into a nice conic program, say, SDP one, on the other.
\par
The goal of what follows is to introduce the notion of ``conic representations'' of problems with convex structure, specifically, convex-concave saddle point problems and monotone VI's, via a given family $\cK$ of cones (e.g., direct products of nonnegative rays/Lorentz/Semidefinite cones). As we will see, given such a representation, the problem of interest can be straightforwardly converted to a conic problem on a cone from $\cK$,  and representations in question admit a fully algorithmic calculus, similar to the calculus of conic representations of convex functions/sets. In  Section \ref{saddlepoint} we {develop} conic representations of convex-concave functions, and {then} in Section \ref{monotone} {extend this framework to} representations of monotone vector fields. {Derivations} in Section \ref{monotone} can be viewed as a ``well-structured'' (i.e., conic representation-oriented) version of convex representation of a monotone operator introduced in \cite[Section 7.4]{NNBook}.

\section{Preliminaries}\label{prelim}
{Consider} a family ${\cK}$ of regular (i.e., closed, convex, pointed, and with nonempty interior) cones in Euclidean spaces such that the family
\begin{enumerate}
\item contains nonnegative ray,
\item is closed w.r.t. taking finite direct products of its members, and
\item is closed w.r.t. passing from a cone to its dual.
\end{enumerate}
{From now on we adopt the following convention: {\em we operate with once for ever fixed family $\cK$ of cones satisfying Assumptions 1--3, and whenever in the sequel we consider a cone, say $\bK$, this cone is assumed to belong to $\cK$.} Besides this, $\bK^*$ stands for the cone dual to $\bK$.}
\par
$\bullet$ {Given a cone $\bK\in \cK$ we } call a constraint in variables $x\in\bR^n$ of the form
\begin{equation}\label{eq00}
Ax\leq_{\bK}a
\end{equation}
{\sl $\cK$-conic}; as usual, $a\leq_{\bK}b$ means that $b-a\in\bK$.  Note that a system
$$
A_kx\leq_{\bK_k}a_k,1\leq k\leq K,
$$
of finitely many $\cK$-conic constraints is equivalent to a single $\cK$-conic constraint{\footnote{We follow the ``Matlab convention'' for matrices: $[A, B]$ and $[A; B]$ denote, respectively, the horizontal and vertical concatenations of two matrices of compatible dimensions.}}
$$
[A_1;...;A_K]x\leq_{\bK}[a_1;...;a_K]
$$
where $\bK=\bK_1\times...\times\bK_K$ belongs to ${\cK}$ along with $\bK_k$, $k\leq K$, since ${\cK}$ is closed w.r.t. taking direct products. For similar reasons, augmenting a $\cK$-conic constraint $Ax\leq_{\bK}a$
in variables $x$ by a system $S$ of finitely many linear equalities and nonstrict linear equalities, we get a conic constraint on $x$. Indeed, $S$ can be written equivalently as $Bx\leq b$, so that the augmented constraint reads
$$
[A;B]x\leq _{\bK\times\bR^m_+}[a;b],\,\,m=\dim b,
$$
and the cone $\bK\times\bR^m_+$ belongs to ${\cK}$ since the latter family contains $\bK$, $\bR_+$ and is closed w.r.t. taking finite direct products.
\par\noindent
$\bullet$ We call conic constraint (\ref{eq00}) {\sl essentially strictly feasible}, if $\bK$ can be decomposed as $\bK=\bR^m_+\times\bL$ with regular cone $\bL$, so that (\ref{eq00}) is of the form
$$
[P;Q]x\leq_{\bR^m_+\times\bL}\leq[p;q],
$$
and there exists $\bar{x}$ such that $P\bar{x}\leq p$ and $Q\bar{x}<_{\bL}q$, the latter meaning that $q-Q\bar{x}\in\inter\, \bL.$
\par\noindent
$\bullet$ We say that a set $X\subset\bR^n$ is {\sl $\cK$-representable} if for properly selected $A,B,c$ and $\bK\in\cF$ one has
$$
X=\{x:\exists u: Ax+Bu\leq_{\bK} c\};
$$
whenever this is the case, the $\cK$-conic constraint $Ax+Bu\leq_{\bK} c$ in the right hand side is called {\sl $\cK$-representation} of $X$. We call a function $f:\bR^n\to\bR\cup\{+\infty\}$ {\sl $\cK$-representable}
if the epigraph $\Epi\{f\}=\{(x,t):t\geq f(x)\}$ is $\cK$-representable {and we refer to a} $\cK$-representation of $\Epi\{f\}$ {as the} $\cK$-representation of $f$. In other words, {asserting} that $\cK$-conic constraint $Ax+tb+Bu\leq_{\bK} c$ represents $f$ is the same as {saying} that the following equivalence takes place:
$$
t\geq f(x)\Leftrightarrow \exists u: Ax+tb+Bu\leq_{\bK} c.
$$
Observe that the level sets $\{x:f(x)\leq a\}$, $a\in\bR$, of $\cK$-representable function $f$ are $\cK$-representable:
$$
\left\{t\geq f(x)\Leftrightarrow \exists u: Ax+tb+Bu\leq_{\bK} c\right\}\Rightarrow \{x: f(x)\leq a\}=\{x:\exists u: Ax+By\leq_{\bK} c-ab\}.
$$
{Thus,} given $\cK$-representations $\{x:\exists u:Ax+Bu\leq_{\bK} c\}$  of $X\subset\bR^n$ and $\{t\geq f(x)\Leftrightarrow \exists v: Cx+td+Ev\leq_{\bL}g\}$ of  $f:\bR^n\to\bR\cup\{+\infty\}$, the problem
{$$\min_{x\in X} f(x)\eqno{(\cP)}$$} is equivalent to {the} conic problem
$$\min_{x,t,u,v}\left\{t:[Ax+Bu;Cx+td+Ev]\leq_{\bK\times\bL} [c;g]\right\}\eqno{(\Q)}
$$
on a cone from $\cK$,
where equivalence of $(\P)$  and $(\Q)$ precisely means that $x$ is a feasible solution to $(\P)$  with a finite value of the objective if and only if there exist $t,u,v$ such that $t\geq f(x)$ and
$x,t,u,v$ is feasible for $(\Q)$.
\par
$\cK$-representable functions/sets admit fully algorithmic calculus: all basic convexity-preserving operations with sets and functions, {(e.g.,} taking finite intersections, images/inverse images under affine mappings, and direct products of sets, or taking linear combinations with nonnegative coefficients {and} convexity-preserving compositions of functions{)} as applied to $\cK$-representable operands, produce $\cK$-representable results, with representations of the results readily given by those of the operands; for details, see, e.g., \cite{LMCOBook,LMCOLN}.
{When a solver for conic problems on cones from $\cK$ is available,} this calculus allows one to recognize that in the problem {$(\cP)$} of interest  objective $f$ and  domain $X$ are $\cK$-representable and find the corresponding $\cK$-representations,\footnote{{Note that} calculus rules are fully algorithmic and thus can be implemented by a compiler, the most famous example being the {\tt CVX} software {\cite{cvx} which can}
 handle ``semidefinite representability'' (family $\cK$ comprised of direct products of semidefinite cones), and as a byproduct---conic quadratic representability.} {so} the problem of interest can be converted into a conic problem {$(\cQ)$} on a cone from $\cK$ and solved by the solver at hand.
\par
From now on, we operate with a once for ever fixed family $\cK$ of cones satisfying the above Assumptions 1--3, which allows for the convention as follows: whenever in the sequel a cone arises, this cone, if the otherwise is explicitly stated, belongs to $\cK$. Besides this, $\bK^*$ stands for the cone dual to $\bK$.
\section{Well-structured convex-concave saddle point problems}\label{saddlepoint}

As  far as its paradigm and set of rules  are concerned,``calculus of $\cK$-representability'' for, say, $\cK$ specified as $\cal{SDP}$ (finite direct products of semidefinite cones)   covers  all {``basic''} needs of ``well-structured'' convex minimization. There is, however, an exception--- the case in which the objective in the convex problem of interest
$$
\min_{x\in \X}\overline{\Psi}(x)\eqno{(P)}
$$
is given implicitly:
\begin{equation}\label{sppsppspp}
\overline{\Psi}(x)=\sup_{y\in\Y}\psi(x,y)
\end{equation}
where $\Y$ is convex set and $\psi:\X\times\Y\to\bR$ is convex-concave (i.e., convex in $x\in\X$ and concave in $y\in\Y$) and continuous. Problem $(P)$ with objective given by (\ref{sppsppspp}) is called ``primal problem associated with the convex-concave saddle point problem
$
\min_{x\in\X}\max_{y\in\Y}\psi(x,y)
$,'' and problems of this type do arise in some applications of well-structured convex optimization. We are about to define a {\sl saddle point} version of $\cK$-representability along with the corresponding calculus which allows to convert ``$\cK$-representable convex-concave saddle point problems'' into usual conic problems on cones from $\cK$.
\subsection{Conic representability of convex-concave function---definition}\label{convconcreprdef}
 Let $\X$, $\Y$ be nonempty convex sets given by  $\cK$-representations:
$$
\X=\{x:\exists \xi: Ax+B\xi\leq_{\bK_\X} c\},\,\, \Y=\{y:\exists \eta: Cy+D\eta\leq_{\bK_\Y} e\}\eqno{[\bK_\X\in\cK,\bK_\Y\in\cK]}.
$$
Let us say that a convex-concave continuous function $\psi(x,y):\X\times\Y\to\bR$ is $\cK$-representable on $\X\times\Y$,
if it admits representation of the form
\begin{equation}\label{cKrepr}
\forall (x\in\X,y\in\Y): \psi(x,y)=\inf\limits_{f,t,u}\left\{f^Ty+t: Pf+tp+Qu+Rx \leq_\bK s\right\}
\end{equation}
where  $\bK\in\cK$. We call representation (\ref{cKrepr}) {\sl essentially strictly feasible}, if the conic constraint
$$
Pf+tp+Qu\leq_\bK s-Rx
$$
in variables $f,t,u$ is essentially strictly feasible for every $x\in\X$.
\subsection{Main observation}
Assume that $\Y$ is compact and is given by essentially strictly feasible $\cK$-representation
\begin{equation}\label{eq16}
\Y=\{y:\exists\eta:Cy+D\eta\leq_{\bK_\Y} e\}.
\end{equation} Then
 problem $(P)$ can be processed as follows: for $x\in\X$ we have
$$
\begin{array}{rl}
\overline{\Psi}(x)=&\max\limits_{y\in\Y}\inf\limits_{f,t,u}\left[f^Ty+t:Pf+tp+Qu+Rx \leq_\bK s\right]\\
=&\inf\limits_{f,t,u}\left\{\max\limits_{y\in\Y}[f^Ty+t]:Pf+tp+Qu+Rx \leq_\bK s \right\}\quad
\left[\hbox{\begin{tabular}{l}Sion-Kakutani Theorem; recall\\
that $\Y$ is convex and compact\\
\end{tabular}}\right]\\
=&\inf\limits_{f,t,u}\left\{\max\limits_{y,\eta}\left[f^Ty:Cy+D\eta\leq_{\bK_\Y} e\right]+t:Pf+tp+Qu+Rx\leq_\bK s\right\}\\
=&\inf\limits_{f,t,u}\left[\min\limits_{\lambda}\left[\lambda^Te:C^T\lambda=f,D^T\lambda=0,\lambda\geq_{\bK_\Y^*}0\right]+t:Pf+tp+Qu+Rx\leq_\bK s\right\}\\
&\multicolumn{1}{r}{\left[\hbox{\begin{tabular}{l}by strong conic duality, see, e.g., \cite[Theorem 1.4.4]{LMCOLN};\\
recall that (\ref{eq16}) is essentially strictly feasible,
\end{tabular}}\right]}\\
\end{array}
$$
so that the problem of interest
$$
\min_{x\in\X} \overline{\Psi}(x)\eqno{(a)}
$$
reduces to the explicit $\cK$-conic problem
$$
\min\limits_{x,\xi,f,t,u,\lambda}\left\{e^T\lambda+t:\begin{array}{l}Pf+tp+Qu+Rx\leq_\bK s,\\
C^T\lambda=f,D^T\lambda=0,\lambda\geq_{\bK_\Y^*}0,\\
Ax+B\xi\leq_{\bK_\X} c\\
\end{array}
\right\}.\eqno{(b)}
$$
{Here,``reduction'' means} that the $x$-component of a feasible solution $\zeta=(x,\xi,f,t,u,\lambda)$ to $(b)$ is a feasible solution to $(a)$ with the value of the objective of the latter problem at $x$ being $\leq$ the value of the objective of $(b)$ at $\zeta$, and the optimal values in $(a)$ and $(b)$ are the same. Thus, as far as building feasible approximate solutions of a prescribed accuracy $\epsilon>0$ in terms of the objective are concerned, problem $(a)$ reduces to the explicit conic problem $(b)$. Note, however, that $(a)$ and $(b)$ are not ``exactly the same''---it may happen that $(a)$ is solvable while $(b)$ is not so. ``For all practical purposes,'' this subtle difference is of no importance since in actual computations exactly optimal solutions usually are not reachable anyway.
\paragraph{Discussion.} Note that for continuous convex-concave function $\psi:\X\times\Y\to\bR$ the set
$$
\Z=\{[f;t;x]: x\in \X,f^Ty+t\geq {\psi}(x,y)\,\forall y\in \Y\}
$$
clearly is convex, and by the standard Fenchel duality we have
\begin{equation}\label{eq444}
\forall (x\in\X,y\in\Y): \psi(x,y)=\inf\limits_{f,t}\left[f^Ty+t:[f;t;x]\in \Z\right].
\end{equation}
$\cK$-representability  of $\psi$ on $\X\times\Y$ means that (\ref{eq444}) is preserved when replacing the set $\Z$ with its properly selected $\cK$-representable subset. Given that $\Z$ is convex, this assumption seems to be not too restrictive; taken together with $\cK$-representability of $\X$ and $\Y$, it can be treated as the definition of $\cK$-representability of the convex-concave function $\psi$. The above derivation shows that convex-concave saddle point problem with $\cK$-representable domain and cost function (more precisely, the primal minimization problem $(P)$ induced by  this saddle point problem) can be represented in explicit $\cK$-conic form, at least when the $\cK$-representations of the cost and of (compact) $\Y$ are essentially strictly feasible.
\par
Note also that if $\X$ and $\Y$ are convex sets and a function $\psi(x,y):\X\times\Y\to\bR$ admits representation (\ref{cKrepr}), then $\psi$ automatically is convex in $x\in\X$ and concave in $y\in\Y$.
\subsection{Symmetry} Assume that representation (\ref{cKrepr}) is essentially strictly feasible. Then for all $x\in\X,y\in\Y$ we have
by conic duality
$$
\begin{array}{rcl}
\psi(x,y)&=&\inf\limits_{f,t,u}\left\{f^Ty+t: Pf+tp+Qu+Rx \leq_\bK s\right\}\\
&=&\sup\limits_{\overline{u}\in\bK^*}\left\{\overline{u}^T[Rx-s]:P^T\overline{u}+y=0,p^T\overline{u}+1=0,Q^T\overline{u}=0\right\},\\
\end{array}
$$
whence, setting
$$
\overline{\X}=\Y,\overline{\Y}=\X,\overline{x}= y,\overline{y}= x,\overline{\psi}(\overline{x},\overline{y})=-\psi(\overline{y},\overline{x})=-\psi(x,y),
$$
we have
$$
\begin{array}{rcl}\multicolumn{3}{l}{(\forall\overline{x}\in\overline{\X},\overline{y}\in\overline{\Y}):}\\
\overline{\psi}(\overline{x},\overline{y})&=&-\psi(x,y)
=\inf\limits_{\overline{u}\in\bK^*}\left\{-\overline{u}^T[Rx-s]:P^T\overline{u}+y=0,p^T\overline{u}+1=0,Q^T\overline{u}=0\right\}\\
&=&\inf\limits_{\overline{f},\overline{t},\overline{u}}\bigg\{
\overline{f}^T\overline{y}+\overline{t}:\underbrace{\left[\begin{array}{l}\overline{f}=-R^T\overline{u},\overline{t}=s^T\overline{u},
Q^T\overline{u}=0,\\
p^T\overline{u}+1=0,P^T\overline{u}+\overline{x}=0,\overline{u}\in\bK^*\\
\end{array}
\right]
}_{\Leftrightarrow \overline{P}\,\overline{f}+\overline{t}\overline{p}+\overline{Q}\overline{u}+\overline{R}\overline{x}\leq_{\overline{\bK}}\overline{s}}\bigg\}\\
\end{array}
$$
with $\overline{\bK}\in\cK$. We see that a (essentially strictly feasible) $\cK$-representation  of convex-concave function $\psi$ on $\X\times\Y$ induces straightforwardly a $\cK$-representation of the ``symmetric entity''---the convex-concave function $\overline{\psi}(y,x)=-\psi(x,y)$ on $\Y\times\X$, with immediate consequences for converting the optimization problem
$$
\sup\limits_{y\in \Y}\left[\underline{\Psi}(y):=\inf\limits_{x\in\X}\psi(x,y)\right]\eqno{(D)}
$$
into the standard conic form.

\subsection{Calculus of conic representations of convex-concave functions}
Representations of the form (\ref{cKrepr}) admit a calculus.
\subsubsection{Raw materials}\label{rawmat}
Raw materials for the calculus are given by \begin{enumerate}
\item Functions $\psi(x,y)=a(x)$, where $a(x)$, $\Dom\, a\supset \X$, is $\cK$-representable:
$$
t\geq a(x) \Leftrightarrow\exists u: \overline{R}x+ t\overline{p}+\overline{Q}u\leq_{\overline{\bK}} \overline{s}
$$
In this case
$$
\psi(x,y)=\inf\limits_{f,t,u}\bigg\{f^Ty+t:\underbrace{f=0,\;\overline{R}x+t\overline{p}+\overline{Q}u\leq_{\overline{\bK}} \overline{s}}_{\Leftrightarrow Pf+tp+Qu+Rx\leq_{\bK} s\hbox{\tiny\ with $\bK\in\cK$}}\bigg\}.
$$
\item Functions $\psi(x,y)=-b(y)$, where $b(y)$, $\Dom\, b\supset \Y$, is $\cK$-representable:
$$
t\geq b(y) \Leftrightarrow\exists u: \overline{R}y+ t\overline{p}+\overline{Q}u\leq_{\overline{\bK}} \overline{s}
$$
with essentially strictly feasible $\cK$-representation. In this case
$$
\begin{array}{rcl}
\psi(x,y)&=&-b(y)=-\inf\limits_{t,u}\left\{t:\overline{R}y+ t\overline{p}+\overline{Q}u\leq_{\overline{\bK}} \overline{s}\right\}\\
&=&-\sup\limits_{u\in\overline{\bK}^*}\left\{-[\overline{R}^Tu]^Ty-s^Tu:-u^T\overline{p}=1,\overline{Q}^Tu=0\right\}\hbox{\ [by conic duality]}\\
&=&\inf\limits_{f,t,u}\bigg\{f^Ty+t: \underbrace{f=R^Tu,t=s^Tu,\overline{p}^Tu+1=0,\overline{Q}^Tu=0,u\geq_{\overline{\bK}^*}0}_{\Leftrightarrow Pf+tp+Qu\leq_{\bK}s\hbox{\tiny\ with $\bK\in\cK$}}\bigg\}.
\end{array}
$$
\item Bilinear functions:
$$
\psi(x,y)\equiv a^Tx+b^Ty+x^TAy+c\Rightarrow \psi(x,y)=\min\limits_{f,t}\bigg\{f^Ty+t:\underbrace{f=A^Tx+b,t=a^Tx+c}_{\Leftrightarrow Pf+tp+Rx \leq s}\bigg\}.
$$
\item\label{GBLF} ``Generalized bilinear functions.'' Let {$\bU\in\cK$ and $E$ be the embedding Euclidean space of $\bU$}.
\begin{enumerate}
\item Let $\overline{\X}$ be a nonempty $\cK$-representable set, and let continuous mapping  $F(x):{\overline\cX}\to E$ possess $\cK$-representable $\bU$-epigraph\footnote{{This implies}, in particular that the $\bU$-epigraph of $F$ is convex, or, which is the same, that $F$ is {\sl $\bU$-convex:}
    $$
    \forall (x';x''\in\overline{\X},\lambda\in[0,1]): F(\lambda x'+(1-\lambda)x'')\leq_{\bU}\lambda F(x')+(1-\lambda) F(x'').
    $$}
\def\Epi{\hbox{\rm Epi}}
$$
\Epi_{\bU}F:=\{(x,z)\in\overline{\X}\times E,z\geq_{\bU} F(x)\}=\{(x,z):\exists u: \overline{R}x+\overline{S}z+\overline{T}u\leq_{\overline{\bK}}\overline{s}\}.
$$
Then the function
$$
\overline{\psi}(x,y)=y^TF(x):\overline{\X}\times\bU^*\to\bR
$$
is $\cK$-representable on $\overline{\X}\times\bU^*$:
$$
\begin{array}{rcl}
\multicolumn{3}{l}{\forall (x\in\overline{\X},y\in\bU^*):}\\
\overline{\psi}(x,y)&=&y^TF(x)=\inf\limits_{f}\left\{f^Ty: f\geq_{\bU}F(x)\right\}\\
&=&\begin{array}{c}
\inf\limits_{f,u}\left\{f^Ty: \overline{R}x+\overline{S}z+\overline{T}u\leq_{\overline{\bK}} \overline{s}\right\}.
\end{array}
\end{array}
$$
\item\label{todayitem} Let  $\Y$ be a nonempty $\cK$-representable set, and let continuous mapping  $G(y):\cY\to E$ possess $\cK$-representable $\bU^*$-hypograph,
\def\Hypo{\hbox{\rm Hypo}}
$$
\begin{array}{rcl}
\Hypo_{\bU^*}G&:=&\{(y,w)\in\Y\times E: w\leq_{\bU^*} G(y)\}\\
&=&\{(y,w):\exists u: \underline{R}y+\underline{S}w+\underline{Q}u\geq_{\underline{\bK}}\underline{s}\}\qquad [\underline{\bK}\in\cK],
\end{array}
$$
the representation being essentially strictly feasible.
Then the function
$$
\underline{\psi}(\underline{x},y)=\underline{x}^TG(y):
\bU\times\Y\to\bR$$ is $\cK$-representable on $\bU\times\Y$:
\begin{equation}\label{today2}
\begin{array}{rcl}
\multicolumn{3}{l}{\forall (\underline{x}\in\bU,y\in\Y):}\\
\underline{\psi}(\underline{x},y)&=&\underline{x}^TG(y)=\sup\limits_{w}\left\{\underline{x}^Tw: w\leq_{\bU^*}G(y)\right\}\hbox{\ [due to $\underline{x}\in\bU$]}\\
&=&\sup\limits_{w,u}\left\{\underline{x}^Tw: \underline{R}y+\underline{S}w+\underline{Q}u\geq_{\underline{\bK}} \underline{s}\right\}\\
&=&\inf_{\lambda}\left\{[\underline{s}-\underline{R}y]^T\lambda:
\underline{S}^T\lambda=\underline{x},\underline{Q}^T\lambda=0,\lambda\in\underline{\bK}^*
\right\}\hbox{\ [by conic duality]}\\
&=&\inf_{f,t,u=[\lambda;w]}\bigg\{f^Ty+t:\underbrace{\left[\begin{array}{l}
f+\underline{R}^T\lambda=0,\underline{s}^T\lambda=t,\\
\underline{S}^T\lambda=\underline{x},\underline{Q}^T\lambda=0,\lambda\in\underline{\bK}^*\\
\end{array}\right]}_{\Leftrightarrow Pf+tp+Qu+R\underline{x}\leq_{\bK} s}\bigg\}\quad[\bK\in\cK].
\end{array}
\end{equation}
\item Let $\Y$ and $G(\cdot)$ be as in item \ref{todayitem}) {with $G(\Y)\subset\bU^*$}, let $\X$ be a nonempty $\cK$-representable set, and let $$F(x):\X\to\bU$$ be continuous $\bU$-convex mapping with $\cK$-representable $\bU$-epigraph:
\begin{equation}\label{today1}
\begin{array}{rcl}
\Epi_\bU F&:=&\{(x,z): x\in\cX,z\geq_{\bU}F(x)\}\\
&=&\{x,z:\exists v: \widehat{R}x+\widehat{S}z+\widehat{Q}v\leq_{\widehat{\bK}}\widehat{s}\}\qquad[\widehat{\bK}\in\cK]\\
\end{array}
\end{equation}
Then the function
$$\psi(x,y)=F^T(x)G(y):\X\times\Y\to \bR$$
is continuous convex-concave and admits $\cK$-representation as follows:
$$
\begin{array}{rcl}
\multicolumn{3}{l}{\forall ({x}\in\X,y\in\Y):}\\
\psi(x,y)&=&F^T(x)G(y)\\
&=&\inf\limits_{z}\left\{z^TG(y):z\geq_\bU F(x)\right\}\hbox{\ [since $G(y)\in\bU^*$]}\\
&=&\inf\limits_{z,v}\left\{z^TG(y):\widehat{R}x+\widehat{S}z+\widehat{Q}v\leq_{\widehat{\bK}}\widehat{s}\right\}\hbox{\ [by (\ref{today1})]}\\
&=&\inf\limits_{z,v}\left\{\inf\limits_{f,t}\left\{f^Ty+t:Pf+tp+Qu+Rz\leq_{\bK} s\right\}:\widehat{R}x+\widehat{S}z+\widehat{Q}v\leq_{\widehat{\bK}}\widehat{s}\right\}
\end{array}
$$
{due to (\ref{today2})---note that on the domain on which $\inf_{z,v}$ is taken we have $z\geq_{\bU} F(x)\in\bU$, making (\ref{today2}) applicable. We conclude that}
\[
\psi(x,y)=
\inf\limits_{f,t,u=[z;v]}\bigg\{f^Ty+t:\underbrace{\left[\begin{array}{l}Pf+tp+Qu+Rz\leq_{\bK} s,\\
\widehat{R}x+\widehat{S}z+\widehat{Q}v\leq_{\widehat{\bK}}\widehat{s}\\
\end{array}\right]}_{\Leftrightarrow \widetilde{P}f+t\widetilde{p}+\widetilde{Q}u+\widetilde{R}x\leq_{\widetilde{\bK}}\widetilde{c}}\bigg\}\\
\quad[\widetilde{\bK}\in\cK].
\]\end{enumerate}
\end{enumerate}
\subsubsection{Basic calculus rules}
Basic calculus rules are as follows.
\begin{enumerate}
\item{[Direct summation]} {Let} $\theta_i>0$, $i\leq I$, and {let}
$$
\begin{array}{c}
\multicolumn{1}{l}{\forall (x^i\in\X^i,y^i\in\Y^i,i\leq I):}\\
 \psi_i(x^i,y^i)=\inf\limits_{f_i,t_i,u_i}\left\{f_i^Ty^i+t_i: P_if_i+t_ip_i+Q_iu_i+R_ix^i \leq_{\bK_i} s_i\right\}\\
 \end{array}.
$$
Then
$$
\begin{array}{rcl}
\multicolumn{3}{l}{\forall \big(x=[x^1;...;x^I]\in\X=\X_1\times...\times \X_I,y=[y^1;...;y^I]\in\Y=\Y_1\times...\times \Y_I\big):}\\
\psi(x,y)&:=&\sum\limits_i\theta_i\psi_i(x^i,y^i)\\
&=&\inf\limits_{f,t,u=\{f_i,t_i,u_i,i\leq I\}}\bigg\{f^Ty+t:\underbrace{\begin{array}{l}f=[\theta_1f_1;...;\theta_I f_I],t=\sum_i\theta_it_i,\\
P_if_i+t_ip_i+Q_iu_i+R_ix^i\leq_{\bK_i} s_i,\,1\leq i\leq I\\
\end{array}}_{\Leftrightarrow Pf+tp+Qu+Rx\leq_{\bK} s\hbox{\tiny\ with $\bK\in\cK$}}\bigg\} \\
\end{array}
$$
\item{[Affine substitution of variables]} {Let}
$$
\begin{array}{c}
\multicolumn{1}{l}{
\forall (\xi\in\cX^+,\eta\in\Y^+):}\\
{\psi_+}(\xi,\eta)=\inf\limits_{f_+,t_+,u_+}\left\{f_+^T\eta+t_+: P_+f_++t_+p_++Q_+u_++R_+\xi \leq_{\bK_+} s_+\right\}\\
\end{array},
$$
and
$$
x\mapsto Ax+b:\X\to \X^+,y\mapsto By+c:\Y\to \Y^+.
$$
Then
{$$
\begin{array}{rcl}
\multicolumn{3}{l}{\forall (x\in\X,y\in\Y):}\\
\psi(x,y)&:=&{\psi_+}(Ax+b,By+c)\\
&=&\inf\limits_{f_+,t_+,u_+}\left\{f_+^T(By+c)+t_+:P_+f_++t_+p_++Q_+u_++R_+[Ax+b]\leq_{\bK_+}s_+\right\}\\
&=&\inf\limits_{f,t,u=[f_+;t_+;u_+]}\bigg\{f^Ty+t:\underbrace{\left[\begin{array}{l}f=B^Tf_+,t=t_++f_+^Tc\\
P_+f_++t_+p_++Q_+u_++R_+Ax\leq_{\bK_+}s_+-R_+b\\
\end{array}\right]}_{\Leftrightarrow Pf+tp+Qu+Rx\leq_\bK s\hbox{ \small with $\bK\in\cK$}}\bigg\}.
\end{array}
$$}\noindent
\item{[Taking conic combinations]} {This rule}, evident by itself, is a combination of the two preceding rules:\\
{Let} $\theta_i>0$ and $\psi_i(x,y):\X\times\Y\to\bR$, $i\leq I$, {be} such that
$$
\begin{array}{c}
\multicolumn{1}{l}{\forall (x\in\X,y\in \Y):}\\
 \psi_i(x,y)=\inf\limits_{f_i,t_i,u_i}\left\{f_i^Ty^i+t_i: P_if_i+t_ip_i+Q_iu_i+R_ix^i \leq_{\bK_i} s_i\right\}\\
 \end{array}.
$$
Then
$$
\begin{array}{rcl}
\multicolumn{3}{l}{\forall (x\in\X,y\in\Y):}\\
\psi(x,y)&:=&\sum\limits_i\theta_i\psi_i(x,y)\\
&=&\inf\limits_{f,t,u=\{f_i,t_i,u_i,i\leq I\}}\bigg\{f^Ty+t:\underbrace{\left[\begin{array}{l}f=\sum_i\theta_if_i,t=\sum_i\theta_it_i,\\
P_if_i+t_ip_i+Q_iu_i+R_ix^i\leq_{\bK_i} s_i,\,1\leq i\leq I\\
\end{array}\right]}_{\Leftrightarrow Pf+tp+Qu+Rx\leq_{\bK} s\hbox{\small\, with $\bK\in\cK$}}\bigg\}.
\end{array}
$$
\item{[Projective transformation in $x$-variable]}  {Let}
$$
\forall (x\in\X,y\in \Y): \psi(x,y)=\inf\limits_{f,t,u}\left\{f^Ty+t: Pf+tp+Qu+R{x} \leq_{\bK} s\right\}.
$$
Then
$$
\begin{array}{c}
\multicolumn{1}{l}{(\forall (\alpha,x):\alpha>0,\alpha^{-1}x\in\X,\forall y\in\Y)}\\
\overline{\psi}((\alpha,x),y):=\alpha\psi(\alpha^{-1}x,y)
=\inf\limits_{f,t,u}\bigg\{f^Ty+t:Pf+tp+Qu+Rx-\alpha s\leq_{\bK}0\bigg\}.\\
\end{array}
$$
\item{[Superposition in $x$-variable]} Let $\X$, $\Y$ be $\cK$-representable,  $\overline{\X}$ be a $\cK$-representable subset of some $\bR^n$, and let
{$\cK\ni\bU$ be a cone in $\bR^n$}. {Furthermore, assume that}
$$
\overline{\psi}(\overline{x},y):\overline{\X}\times\Y\to\bR
$$
{is a} continuous convex-concave function which is $\bU$-nondecreasing in $\overline{x}\in\overline{\X}${, i.e.}
$$\forall (y\in\Y,\overline{x}',\overline{x}''\in\overline{\X}:\overline{x}'\leq_{\bU}\overline{x}''):
\overline{\psi}(\overline{x}'',y)\geq\overline{\psi}(\overline{x}',y),
$$
 and admits $\cK$-representation on $\overline{\X}\times\Y$:
$$
\forall(\overline{x}\in\overline{\X},y\in\Y):\,\overline{\psi}(\overline{x},y)=\inf\limits_{{f},{t},\overline{u}}
\left\{{f}^Ty+{t}:\overline{P}{f}+{t}\overline{p}+\overline{Q}\overline{u}+\overline{R}\overline{x}
\leq_{\overline{\bK}}\overline{s}\right\}.
$$
Let also
$$
x\mapsto X(x):\X\mapsto\overline{\X}
$$
be a $\bU$-convex mapping such that the intersection of  the $\bU$-epigraph of the mapping with $\X\times\overline{\X}$ admits $\cK$-representation:
$$
\{(x,\overline{x}):x\in \X,\overline{x}\in\overline{\X},\overline{x}\geq_{\bU} X(x)\}=\{(x,\overline{x}):\exists v: Ax+B\overline{x}+Cv\leq_{\widehat{\bK}} d\}.
$$
Then the function $\psi(x,y)=\overline{\psi}(X(x),y)$ admits $\cK$-representation on $\X\times\Y$:
$$
\begin{array}{rcl}
\multicolumn{3}{l}{\forall (x\in\X,y\in\Y):}\\
\psi(x,y)&=&\overline{\psi}(X(x),y)=\inf\limits_{\overline{x}}\{\overline{\psi}(\overline{x},y):\overline{x}\in\overline{\X}\ \&\ \overline{x}\geq_{\bU} X(x)\}\\
&=&\inf\limits_{{f},{t},\overline{u},\overline{x}}\left\{{f}^Ty+{t}:
\begin{array}{l}\overline{x}\in\overline{\X},\overline{x}\geq_{\bU} X(x)\\
\overline{P}{f}+{t}\overline{p}+\overline{Q}\overline{u}+\overline{R}\overline{x}\leq_{\overline{\bK}}\overline{s}\\
\end{array}\right\}\\
&=&\inf\limits_{{f},{t},\overline{u},\overline{x},v}\left\{{f}^Ty+{t}:
\begin{array}{l}Ax+B\overline{x}+Cv\leq_{\widehat{\bK}}d\\
\overline{P}{f}+{t}\overline{p}+\overline{Q}\overline{u}+\overline{R}\overline{x}\leq_{\overline{\bK}}\overline{s}\\
\end{array}\right\}\\
&=&\inf\limits_{{f},{t},u=[\overline{u},\overline{x},v]}\bigg\{{f}^Ty+{t}:\underbrace{
\left[\begin{array}{l}Ax+B\overline{x}+Cv\leq_{\widehat{\bK}}d\\
\overline{P}{f}+{t}\overline{p}+\overline{Q}\overline{u}+\overline{R}\overline{x}\leq_{\overline{\bK}}\overline{s}\\
\end{array}\right]}_{\Leftrightarrow Pf+tp+Qu+Rx\leq_{\bK}s\hbox{\small\ with $\bK\in\cK$}}\bigg\}.
\end{array}
$$
\item{[Partial maximization]} Let
$$
\begin{array}{l}
\forall (x\in\X,y=[w;z]\in \Y): \\
\psi(x,[w;z])=\inf\limits_{[g;h],\tau,u}\left\{g^Tw+h^Tz+\tau: Gg+Hh+\tau p+Qu+R{x} \leq_{\bK} s\right\}\qquad[\bK\in\cK],\end{array}
$$
and let $\Y$ be compact and given by $\cK$-representation:
$$
\Y=\{[w;z]:\exists v: Aw+Bz+Cv\leq_{\bL} r\}
$$
such that the conic constraint $Bz+Cv\leq_{\bL}r-Aw$ in variables $z,v$ is essentially strictly feasible for every $w\in{\cal W}=\{w:\exists z: [w;z]\in\Y\}$. Then the function
\[\overline{\psi}(x;w):=\max\limits_z\left\{\psi(x,[w;z]):[w;z]\in \Y\right\}:\X\times{\cal W}\to\bR\] is $\cK$-representable { provided it is continuous:\footnote{Representation to follow holds true without the latter assumption; we make it to stay consistent with our general definition of representability, where the convex-concave function in question is assumed to be continuous.}}
$$
\begin{array}{l}
{\forall (x\in\X,w\in{\cal W}):}\\
\overline{\psi}(x;w)=\max\limits_z\Big\{\inf\limits_{[g;h],\tau,u}\left[g^Tw+h^Tz+\tau: Gg+Hh+\tau p+Qu+R{x} \leq_{\bK} s\right]:[w;z]\in\Y\Big\}\\
=\inf\limits_{[g;h],\tau,u}\Big\{\max\limits_z\left[g^Tw+h^Tz+\tau: [w;z]\in\Y\right]:Gg+Hh+\tau p+Qu+R{x} \leq_{\bK} s\Big\}\\
\qquad\hbox{[by the Sion-Kakutani Theorem; note that for $w\in{\cal W}$}\\
\qquad\hbox{the set $\{z:[w;z]\in\Y\}$ is nonempty, convex and compact]}\\
=\inf\limits_{[g;h],\tau,u}\Big\{\max\limits_{z,v}\left[g^Tw+h^Tz+\tau: Bz+Cv\leq_{\bL} r-Aw\right]:\;
{Gg+Hh+\tau p+Qu+R{x} \leq_{\bK} s\Big\}}\\
=\inf\limits_{[g;h],\tau,u}\Big\{\min\limits_{\xi}\left[g^Tw+(r-Aw)^T\xi+\tau: B^T\xi=h,C^T\xi=0,\xi\geq_{\bL^*}0\right]:\\
\qquad\qquad\qquad\qquad\qquad\qquad\qquad\qquad{Gg+Hh+\tau p+Qu+R{x} \leq_{\bK} s\Big\}}\qquad\hbox{\ [by conic duality]}\\
=\inf\limits_{f,t,\overline{u}=[g,h,\tau,u,\xi]}\bigg\{f^Tw+t:\underbrace{\left[\begin{array}{l}f=g-A^T\xi,t=r^T\xi+\tau,\;B^T\xi=h,\;c^T\xi=0\\
\xi\geq_{\bL^*}0,Gg+Hh+\tau p+Qu+R{x} \leq_{\bK} s
\end{array}\right]}_{\Leftrightarrow
\overline{P}f+t\overline{p}+\overline{Q}\overline{u}+\overline{R}x\leq_{\overline{\bK}}\overline{s}\ \hbox{with $\overline{\bK}\in\cK$}}\bigg\}.
\end{array}
$$
Note that the last three rules combine with symmetry to induce ``symmetric'' rules on perspective transformation and superposition in $y$-variable and partial minimization in $x$-variable.
\item{[Taking Fenchel conjugate]} Let $\X\subset\bR^n$, $\Y\subset\bR^m$ be nonempty convex compact sets given by essentially strictly feasible $\cK$-representations
$$
\X=\{x:\exists \xi: Ax+B\xi\leq_{\bK_\X} c\},\,\Y=\{y:\exists \eta:Cy+D\eta\leq_{\bK_\Y} e\},
$$
and assume that the conic constraint
$$
D^T\lambda=0\ \&\ \lambda\geq_{\bK_\Y^*}0
$$
is essentially strictly feasible (this definitely is the case when $\bK_\Y$ is polyhedral). Let, next,
$\psi(x,y):\X\times\Y\to\bR$ be continuous convex-concave function given by essentially strictly feasible $\cK$-representation:
$$
\psi(x,y)=\inf\limits_{f,t,u}\left\{f^Ty+t:Pf+t\Pi+Qu+Rx\leq_{\bK}s\right\}.
$$
Consider the {\sl Fenchel conjugate} of $\psi$---the function
$$
\psi_*(p,q)=\max\limits_{x\in\X}\min\limits_{y\in\Y}\left[p^Tx+q^Ty-\psi(x,y)\right]:\bR^n\times\bR^m\to\bR.
$$
We claim that $\psi_*$ is a continuous convex-concave $\cK$-representable function with $\cK$-representation readily given by the $\cK$-representations of $\X$, $\Y$, $\psi$
\par\noindent
The fact that $\psi_*$ is well defined and continuous is readily given by compactness of $\X$, $\Y$ and continuity of $\psi$. These properties of the data imply by the Sion-Kakutani Theorem  that
$$
\psi_*(p,q)=\min\limits_{y\in\Y}\max\limits_{x\in\X}\left[p^Tx+q^Ty-\psi(x,y)\right].
$$
From the initial $\max\limits_x\min\limits_y$ definition of  $\psi_*$ it follows that
$$
\psi_*(p,q)=\max\limits_{x\in\X}\bigg[p^Tx+\min\limits_{y\in\Y}[q^Ty-\psi(x,y)]\bigg]
$$
is the pointwise maximum of a family of affine functions of $p$ and thus is convex in $p$. From the $\min\limits_y\max\limits_x$ representation of $\psi_*$ it follows that
$$
\psi_*(p,q)=\min\limits_{y\in\cY}\big[q^Ty+\max\limits_{x\in \X}[p^Tx-\psi(x,y)]\big]
$$
is the pointwise minimum of a family of affine functions of $q$ and thus is concave in $q$.\par
It remains to build $\cK$-representation of $\psi_*$. We have
\begin{eqnarray}
\lefteqn{\psi_*(p,q)=\max\limits_{x\in\X}\min\limits_{y\in\Y}\left[p^Tx+q^Ty-\psi(x,y)\right]}
\nonumber\\ \nonumber&=&
\max\limits_{x\in\X}\bigg[p^Tx+\min\limits_{y\in\Y}\sup\limits_{f,t,u}\left\{[q-f]^Ty-t:Pf+t\Pi+Qu+Rx\leq_{\bK} s\right\}\bigg]\\ \nonumber
&=&\max\limits_{x\in\X}\bigg[p^Tx+\sup\limits_{f,t,u}\bigg\{\min\limits_{y\in\Y}\{[q-f]^Ty\}-t:Pf+t\Pi+Qu+Rx\leq_{\bK} s\bigg\}\bigg]\\ \nonumber
&&\hbox{\small[by the Sion-Kakutani Theorem; recall that $\Y$ is convex and compact} \\ \nonumber&&\hbox{and $[q-f]^Ty-t$ is {concave in} $f,t,u$ and {convex} in $y$]}\\ \nonumber
&=&\max\limits_{x\in\X}\bigg[p^Tx+\sup\limits_{f,t,u}\bigg\{\sup\limits_\lambda\left\{-e^T\lambda:f-C^T\lambda=q,D^T\lambda=0,\lambda\geq_{\bK_\Y^*}0\right\}-t:
\\ \nonumber
&&\qquad\qquad\qquad\qquad\qquad\qquad\qquad\qquad\qquad{Pf+t\Pi+Qu+Rx\leq_{\bK} s\bigg\}\bigg]}\\ \nonumber
&&\hbox{\small[by Conic Duality; note that $\Y$ is given by essentially strictly feasible $\cK$-representation]}\\ \nonumber
&=&\sup\limits_{x\in\X,f,t,u,\lambda}\left\{p^Tx-t-e^T\lambda:\begin{array}{l}D^T\lambda=0,\;
f-C^T\lambda=q,
\lambda\geq_{\bK_\Y^*}0\\
Pf+t\Pi+Qu+Rx\leq_{\bK} s
\end{array}\right\}\\
&=&\sup\limits_{x,\xi,f,t,u,\lambda}\left\{p^Tx-t-e^T\lambda:
\begin{array}{ll}D^T\lambda=0&(a)\\
f-C^T\lambda=q&(b)\\
\lambda\geq_{\bK_\Y^*}0&(c)\\
Pf+t\Pi+Qu+Rx\leq_{\bK} s&(d)\\
Ax+B\xi\leq_{\bK_\X} c&(e)\\
\end{array}\right\}
\label{theend?}
\end{eqnarray}
Now, the $\cK$-representation of $\cX$ is essentially strictly feasible, so that $(e)$ admits an essentially strictly feasible solution $\bar{x},\bar{\xi}$; observe that $\bar{x}\in\X$. The $\cK$-representation of $\psi$ is essentially strictly feasible, implying that $\bar{x}$ can be augmented by $\bar{f}$, $\bar{t}$, $\bar{u}$ in such a way that $(\bar{x},\bar{f},\bar{t},\bar{u})$ is an essentially strictly feasible solution to $(d)$. By the origin of constraints $(a)$---$(c)$, their system is feasible for all $f,q$ as a system of constraints on $\lambda$, .
Besides this, by assumption, there exists a representation $\bK_\Y^*=M\times N$ with regular cone $M$ and polyhedral cone $N$ and $\lambda^\prime\in[\inter M]\times N$ such that $D^T\lambda^\prime=0$. Taking into account that the system $(a)$---$(c)$, considered as a system in variable $\lambda$, is solvable for all $f,q$, there exists $\lambda^{\prime\prime}\in\bK_{\Y^*}$ such that $D^T\lambda^{\prime\prime}=0$ and $\bar{f}-C^T\lambda^{\prime\prime}=q+C^T\lambda^\prime$. This implies that $\bar{\lambda}=\lambda^\prime+\lambda^{\prime\prime}$ taken together with $\bar{f}$ is an essentially strictly feasible solution to the system $(a)$---$(c)$ of constraints in variables $\lambda$ and $f$.
The bottom line is that the constraints $(a)$---$(e)$ in variables $x,\xi,f,t,u,\lambda$ form an essentially strictly feasible conic constraint. Besides this,  problem (\ref{theend?}) is bounded  by its origin. By  Conic Duality we conclude from (\ref{theend?}) that
\begin{eqnarray*}
\psi_*(p,q)&=&\min\limits_{\alpha,\beta,\gamma,\delta,\epsilon}\left\{\beta^Tq+s^T\delta+c^T\epsilon:\begin{array}{l}\gamma\geq_{\bK_\Y}0,\delta\geq_{\bK^*}0,\epsilon\geq_{\bK_\X^*}0\\
R^T\delta+A^T\epsilon=p,B^T\epsilon=0,P^T\delta+\beta=0\\
\Pi^T\delta=-1,Q^T\delta=0,D\alpha-C\beta+e=\gamma\\
\end{array}\right\}\\
&=&\min\limits_{\beta,\tau,w=(\alpha,\delta,\epsilon)}\left\{\beta^Tq+\tau:
\begin{array}{l}\delta\geq_{\bK^*}0,\epsilon\geq_{\bK_\X^*}0,\;\tau=s^T\delta+c^T\epsilon\\
R^T\delta+A^T\epsilon=p,B^T\epsilon=0,P^T\delta+\beta=0\\
\Pi^T\delta=-1,Q^T\delta=0,D\alpha-C\beta+e\geq_{\bK_\Y}0\\
\end{array}\right\}
\end{eqnarray*}
which is the desired $\cK$-representation of $\psi_*$. \qed
\end{enumerate}
\subsection{Illustrations} {\bf A.} Our first illustration {is} motivated by {a} statistical application {of saddle point optimization---}near-optimal recovery of linear forms in Discrete observation scheme, see \cite[Section 3.1]{StatOpt}{. Let}
$$
\psi(x,y)=\ln\left(\sum_i\e^{x_i}y_i\right):\X\times\Y\to\bR,
$$
 $\X$ and $\Y$ {be} $\cK$-representable, and {let} $\Y$, $0\not\in\Y$, {be} is a compact subset of the nonnegative orthant. {Because for $z>0$
 \[\ln z=\inf_{u} ze^u-u-1,
\]
 for} $y\geq0$ we clearly have
$$
\begin{array}{c}
\ln\left(\sum_i\e^{x_i}y_i\right)=\inf\limits_{u}\left[\left(\sum_i\e^{x_i}y_i\right)\e^u-u-1\right]=\inf\limits_{f,u}\left[\sum_iy_if_i-u-1:f_i\geq \e^{x_i+u}\right]\\
=\inf\limits_{f,t,u}\left\{f^Tu+t:f_i\geq \e^{x_i+u}\,\forall i\ \&\ t\geq -u-1\right\}
\end{array}
$$
The resulting representation is a $\cK$-representation, provided that the closed w.r.t. taking finite direct products and passing to the dual cone family $\cK$  of regular cones contains $\bR_+$, the exponential cone
$$
\bEx{=}\cl\{[t;s;r]: t\geq s\e^{r/s},s>0\},
$$
and, therefore, its dual cone
$$
\bEx^*=\cl\{[\tau;\sigma;-\rho]:\tau>0,\rho>0,\sigma\geq \rho\ln(\rho/\tau)-\rho\}.
$$
{\bf B.} {Let now}
$$
\psi(x,y)=\left(\sum_{i=1}^n\theta^p_i(x)y_i\right)^{1/p}
$$
where $p>1$, $\theta_i(x)$ are nonnegative $\cK$-representable real-valued functions on $\cK$-representable set $\X$,  and $\Y$ is a $\cK$-representable subset of the nonnegative orthant. In this case, as is easily seen, for all $(x\in\X,y\in\Y)$ it holds
$$
\psi(x,y)=\inf\limits_{[f;t]}\left\{f^Ty+t: t\geq0,f\geq0,t^{{p-1\over p}}f_i^{{1\over p}}\geq \kappa \theta_i(x),i\leq n\right\} \eqno{[\kappa=p^{-1}(p-1)^{{p-1\over p}}]}
$$
which can immediately  be converted into $\cK$-representation, provided $\cK$ contains 3D Lorentz cone
{$\bL^3=\{x\in\bR^3:x_3\geq\sqrt{x_1^2+x_2^2}\}$}
and $p$ is rational, see \cite[Section 3.3]{LMCOBook} or \cite[Section 2.3.5]{LMCOLN}.
\\
{\bf C.} {In our next example,} $\cX\subset\bR^{m\times n}$ and $\Y\subset\bS^m_+$ are nonempty convex sets, and
$$
\psi(x,y)=2\sqrt{\Tr(x^Tyx)}:\X\times\Y\to\bR.
$$
Taking into account that for $a\geq0$ one has $2\sqrt{a}=\inf_{s>0}[a/s+s]$, we have
$$
\begin{array}{rcl}
\multicolumn{3}{l}{\forall (x\in\X,y\in\Y):}\\
\psi(x,y)&=&2\sqrt{\Tr(y[xx^T])}=\inf_{g}\left\{2\sqrt{\Tr(yg)}:g\succeq xx^T\right\}\\
&=&\inf_{f,s}\left\{\Tr(yf)+s:s>0,fs\succeq xx^T\right\}\\&=&
\inf\limits_{f,s}\left\{\Tr(yf)+s:\left[\begin{array}{c|c}f&x\cr\hline x^T&sI_n\cr\end{array}\right]\succeq0\right\}.
\end{array}
$$
The resulting representation is $\cK$-representation, provided that $\cK$ contains semidefinite cones.
\begin{quote}
 This is how {\bf C} works in Robust Markowitz Portfolio Selection {(cf, e.g., \cite{ElOk2003,GI03})}
$$
\min\limits_{x\in\X}\max_{y\in\Y} \left[-r^Tx+2\rho\sqrt{x^Tyx}\right]\eqno{[\rho>0]}
$$
(here $x\in\bR^n$ is the composition of portfolio, $r$ is the vector of expected returns, and $y$ is the uncertain covariance matrix of the returns). Assuming for the sake of definiteness that
$\Y$ is cut off $\bS^n_+$ by the constraints
$$
a_i^Tya_i+b_iy+yb_i^T\preceq p_i,i\leq I, y_-\leq y\leq y_+
$$
(where $\leq$ for matrices acts entrywise) and applying our machinery on the top of the above semidefinite representation of $2\sqrt{x^Tyx}$, the saddle point problem reduces to
$$
\begin{array}{l}
\min\limits_{x,s,\alpha_i,\mu_\pm}\bigg\{-r^Tx+\rho\left[s+\sum_i\Tr(\alpha_ip_i)+\Tr(\mu_+y_+-\mu_-y_-)\right]:\\
\multicolumn{1}{r}{\qquad\qquad\left.\begin{array}{l}\left[\begin{array}{c|c}
\sum_i\left[a_i\alpha_ia_i^T+\alpha_ib_i+b_i^T\alpha_i\right]+\mu_+-\mu_-&x\cr\hline x^T&s\cr\end{array}\right]\succeq0\\
\alpha_i\succeq0,i\leq I,\mu_\pm\geq0,x\in\X
\end{array}\right\}.}
\end{array}
$$
\end{quote}
{\bf D.} In our concluding example, $\cK$ contains the products of semidefinite cones, $\X=\bR^{m\times n}$, $\Y=\bS^n_+$, and
$$\psi(x,y)=\Tr(x^Txy^{1/2}):\X\times \Y\to\bR.$$ This is a ``generalized bilinear function''; in terms of item \ref{GBLF}.c of Section \ref{rawmat}, we have $F(x)=x^Tx$, $G(y)=y^{1/2}$, $\bU=\bU^*=\bS^n_+$, and
$$
\begin{array}{rcl}
\Epi_\bU F&:=&\{(x,z):z\succeq x^Tx\}=\left\{(x,z):\left[\begin{array}{c|c}z&x^T\cr\hline x&I_m\cr\end{array}\right]\succeq0\right\},\\
\mathrm{Hypo}_{\bU^*}G&:=&\{(y,w):y\in\bS^n_+,w\preceq y^{1/2}\}=\left\{(y,w):\exists v:\left[\begin{array}{c|c}y&v\cr\hline v&I_n\cr\end{array}\right]\succeq0,v\succeq 0,w\preceq v\right\}.
\end{array}
$$
With these data, the construction from item \ref{GBLF}.c of Section \ref{rawmat} leads straightforwardly to the following semidefinite representation of $\psi$:
{\small$$
\forall (x\in\bR^{m\times n},y\in\bS^n_+): \psi(x,y):=\Tr(x^Txy^{1/2})=\inf\limits_{f,t,u=(z,\beta,\gamma)}\left\{\Tr(fy)+t:
\begin{array}{l}f\in\bS^n,\beta\in\bR^{n\times n},\gamma\in\bS^n,z\in\bS^n\\
t=\Tr(\gamma), z\preceq \beta+\beta^T\\
\left[\begin{array}{c|c}f&\beta\cr\hline\beta^T&\gamma\cr\end{array}\right]\succeq0,\left[\begin{array}{c|c}z&x^T\cr\hline x&I_m\cr\end{array}\right]\succeq0\\
\end{array}\right\}.
$$}
\section{Well-structured variational inequalities with monotone operators}\label{monotone}
\subsection{Preliminaries}
Let $\X$ be a nonempty convex subset of Euclidean space $E$. A {\sl vector field} $F(x):x\to E$ is called {\sl monotone}, if
$$
\langle F(x)-F(y),x-y\rangle \geq0\;\forall x,y\in\X.
$$
{\sl Variational Inequality} $\VI(F,\X)$ reads
$$
\hbox{Find\ }x_*\in \X: \langle F(y),{x_*}-y\rangle\leq0\;\forall y\in\X.\eqno{\VI(F,\X)}
$$
Solutions which are sought in $\VI(F,\X)$ are called {\sl weak solutions}; they do exist whenever $\X$ is a compact { set and $F$ is monotone on $\X$}. {\em Strong solutions} are points $x_*\in\X$ such that
$$
\langle F(x_*),x-x_*\rangle\geq0\;\forall x\in\X;
$$
from monotonicity of $F$ it follows that every strong solution is a weak one. The opposite is true provided that $F$ is continuous on $\X$.
\par
A natural (in)accuracy measure for a candidate approximate solution to $\VI(F,\X)$ is the {\sl dual gap function}
$$
\epsVI[x|\X]=\sup_{y\in \X} \langle F(y),x-y\rangle;
$$
this function is nonnegative and is zero exactly at weak solutions to $VI(F,\X)$.
\par
Our current goal is to develop {a framework} for converting ``well-structured'' VI's with monotone vector fields into the usual conic problems, thus bringing them within the grasp of Interior-Point polynomial time methods and software like {\tt CVX}. What follows can be seen as a streamlined ``well-structured'' (i.e., conic representation-oriented) extension of what in \cite[Section 7.4]{NNBook} was called a convex representation of a monotone operator.

\subsection{Conic representability of monotone vector fields and monotone VI's in conic form}
\subsubsection{Conic representation of a monotone vector field} \label{cKrMVF}
Let $F(x):\X\to{E}$ be a monotone and continuous vector field on nonempty convex subset $\X$ of Euclidean space $E$. Consider the set
\begin{equation}\label{jan333}
{\cal F}[F,\X]=\{[t;g;x]\in \bR\times E\times E: x\in\X, t-\langle g,y\rangle \geq \langle F(y),x-y\rangle\,\,\forall y\in\X\}
\end{equation}
and let us make several straightforward observations.
\par\noindent
{\bf\ref{cKrMVF}.A.} ${\cal F}[F,\X]$ is a convex set which contains all triples $[\langle{F(x)},x\rangle;F(x);x]$, $x\in\X$; this set is closed provided $\X$ is so. Besides this, $\cF[F,\X]$ is $t$-monotone:
$$[t;g;x]\in\cF[F,\X]\;\mbox{and}\; t'\geq t\Rightarrow [t';g;x]\in\cF[F,\X].
$$
\begin{quote}Indeed, ${\cal F}[F,\X]$ is the intersection of solution set of the system of nonstrict linear constraints
$$
t-\langle g,y\rangle -\langle F(y),x\rangle \geq -\langle F(y),y\rangle,\,y\in\X
$$ in variables $g,t,x$ (this set is closed and convex) with the convex set $\{[t;g;x]:x\in\X\}$ (which is closed if so is $\X$). By monotonicity, for every $x\in \X$ we have
$$
\langle F(x),x-y\rangle \geq \langle F(y),x-y\rangle \,\,\forall y\in\X,
$$
so that $[\langle F(x),x\rangle;F(x);x]\in\cF[F,\X]$; and $t$-monotonicity is evident.
\end{quote}
\noindent
{\bf \ref{cKrMVF}.B.} For $\epsilon\geq0$, let \begin{equation}\label{xstar}
\X_*(\epsilon)=\{[g;t]\in E\times\bR: \sup\limits_{y\in\X}[t-\langle g,y\rangle]\leq\epsilon\},
\end{equation}
 so that $\X_*(\epsilon)$ is a nonempty closed convex set. Then for every $\epsilon\geq0$,
$\epsilon$-solutions to  $\VI(F,\X)$---points  $x\in \X$ such that
$$
\epsVI[x|F]:=\sup\limits_{y\in\X}\langle F(y),x-y\rangle \leq\epsilon
$$
are exactly the points
\begin{equation}\label{jan99}
x:\exists(t,g): [t;g;x]\in\cF[F,\X]\ \;\mbox{and}\;\ [t;g]\in\cX_*(\epsilon);
\end{equation}
\begin{quote}
Indeed, let $x,t,g$ be such that $[t;g;x]\in\cF[F,\cX]$ and $[t;g]\in\X_*(\epsilon)$. From the first inclusion it follows that $x\in\X$ and
$$
t-\langle g,y\rangle \geq \langle F(y),x-y\rangle \,\,\forall y\in\X,
$$
while from the second inclusion it follows that $t-\langle g,y\rangle\leq\epsilon$ for all $y\in \X$. Taken together, these two relations
imply that $x\in \X$ and $\langle F(y),x-y\rangle \leq\epsilon$ for all $y\in \X$,  i.e., $\epsVI(x|F)\leq\epsilon$. Vice versa, if $x\in\X$, $\epsilon\geq0$, and $\epsVI[x|F]\leq\epsilon$, then
the triple $t=\epsilon,g=0,x$ clearly belongs to $\cF[F,\X]$ and $[t;g]\in\X_*(\epsilon)$.
\end{quote}
\paragraph{$\cK$-representation of $(F,\X)$. }
Given a continuous monotone vector field $F:\X\to E$ on a nonempty convex subset $\X$ of Euclidean space $E$, let us call a conic constraint{\footnote{Recall that we have fixed a family $\cK$ of regular cones in Euclidean spaces, and by our standing convention, all cones to be considered below belong to $\cK$.}}
\begin{equation}\label{jan0}
Xx+Gg+tT+Uu\leq_{\bK} a
\end{equation}
in variables $t\in\bR$, $g\in E$, $x\in E$, and $u\in\bR^k$
{\sl a $\cK$-representation} of $(F,\X)$, if the set
\begin{equation}\label{jan100}
\T:=\{[t;g;x]:\exists u: Xx+Gg+tT+Uu\leq_{\bK} a\}
\end{equation}
possesses the following two properties:
\begin{itemize}
\item[(i)] $\T$ is contained in $\cF[F,\X]$ and, along with the latter set, is ``$t$-monotone:''
$$
[t;g;x]\in \T\ \;\mbox{and}\;\ t'\geq t\Rightarrow [t';g;x]\in\cT
$$
(note that when $\T\neq\emptyset$, $t$-monotonicity is equivalent to ${T}\leq_{\bK}0$);
\end{itemize}
and
\begin{itemize}
\item[(ii)] $\T$ contains the set
$$
\{[\langle F(x),x\rangle;F(x);x],x\in\X\}\eqno{[\subset \cF[F,\X]]}
$$
\end{itemize}
If the set (\ref{jan100}) satisfies (i) and the relaxed version of (ii), specifically,
\begin{itemize}
\item[(ii')]
$\T$ contains the set
$$
\{[t;F(x);x]: t>\langle F(x),x\rangle,x\in\X\}\eqno{[\subset \cF[F,\X]]}
$$
\end{itemize}
we say that (\ref{jan0})  is an {\sl almost $\cK$-representation} of $(F,\X)$.
\par
Let us make an immediate observation:
{\begin{remark}\label{rem1} In the situation described in the beginning of this section, let $\X$ be closed, and let $\Y$ be a convex set such that
$$
\Conv\{[x^TF(x);F(x);x],x\in\X\}\subset \Y\subset \cl\Conv\{[x^TF(x);F(x);x],x\in\X\}
$$
and let $\Y_+=\Y+\bR_+\times\{0_E\}\times\{0_E\}$, where $0_E$ is the origin in $E$.
Then every $\cK$-representation of $\Y_+$ represents $(F,\X)$.
\end{remark}}\noindent
\subsubsection{Conic form of conic-representable monotone VI}
Our main observation is as follows:
\begin{proposition}\label{propreprjan}
Let $\X\subset E$ be nonempty convex compact given by essentially strictly feasible $\cK$-representation:
\begin{equation}\label{jan40}
\X=\{x:\exists v: Ax+Bv\leq_{\bL} b\},\qquad\qquad{[\bL\in\cK],}
\end{equation}
so that $\X_*(\epsilon)$, see {\rm (\ref{xstar})},  by conic quality, admits $\cK$-representation as follows:
\begin{equation}\label{jan41}
\X_*(\epsilon)=\{[t;g]: \exists\lambda:A^T\lambda+g=0,B^T\lambda=0,t+b^T\lambda\leq\epsilon,\lambda\geq_{\bL^*}0\}.
\end{equation}
Let, moreover, $F:\X\to E$ be a continuous monotone vector field, and let {\rm (\ref{jan0})} be an almost $\cK$-representation of $(F,\X)$.
Then for every $\epsilon>0$ the  system of conic constraints
\begin{equation}\label{jan42}
\begin{array}{rcll}
Xx+Gg+tT+Uu&\leq_{\bK}& a&(a)\\
A^T\lambda+g&=&0&(b.1)\\
B^T\lambda&=&0&(b.2)\\
t+b^T\lambda&\leq&\epsilon&(b.3)\\
\lambda&\geq_{\bL^*}&0&(b.4)\\
\end{array}
\end{equation}
in variables $x,g,t,u,$ and $\lambda$ is feasible, and $x$-component of any feasible solution belongs to $\X$ and is an $\epsilon$-solution to $\VI(F,\X)$:
$$
\epsVI[x|F]\leq\epsilon.
$$
Therefore, finding {an} $\epsilon$-solution to $\VI(F,\X)$ reduces to finding a feasible solution to an explicit $\cK$-conic constraint.\par
When {\rm (\ref{jan0})} is a $\cK$-representation of $(F,\X)$, the above conclusion holds true for all $\epsilon\geq0$.
\end{proposition}
\noindent
{\bf Proof.} Let us fix $\epsilon>0$, and let $\bar{x}$ solve $\VI(F,X)$ (a solution exists since $\X$ is compact). Then $\langle F(\bar{x}),\bar{x}-y\rangle\leq0$ for all $y\in \X$ due to continuity of $F$ on $\X$.
Consequently, given $\epsilon>0$, setting $\bar{t}=\langle F(\bar{x}),\bar{x}\rangle$, $\bar{g}=F(\bar{x})$, we have $[\bar{t}+\epsilon;\bar{g};\bar{x}]\in \X_*(\epsilon)$, implying by (\ref{jan41}) that there exists $\bar{\lambda}$ such that $t=\bar{t}+\epsilon$, $g=\bar{g}$, $\lambda=\bar{\lambda}$ solve (\ref{jan42}.$b$.1-4). Besides this, by (ii') there exists $\bar{u}$ such that $x=\bar{x}$, $g=\bar{g}$, $t=\bar{t}+\epsilon$, $u=\bar{u}$ solve (\ref{jan42}.$a$).   Thus, (\ref{jan42}) is feasible.
\par
Next, if $x,g,t,u,\lambda$ is a feasible solution to (\ref{jan42}), then by (\ref{jan42}.$a$) one has $[t;g;x]\in\T$, with $\T$ given by (\ref{jan100}), whence $[t;g;x]\in\cF[F,\X]$, implying that $x\in\X$, see (\ref{jan333}). Besides this, (\ref{jan42}.$b$.1-4) say that
$[t;g]\in\X_*(\epsilon)$.  Thus, $[t;g;x]\in\cF[F,\X]$ and $[t;g]\in \X_*(\epsilon)$, implying by \ref{cKrMVF}.B that $\epsVI[x|F]\leq\epsilon$.
\par
Finally, when (\ref{jan0}) is a $\cK$-representation of $(F,\X)$, the above reasoning works for $\epsilon=0$. \qed
\subsection{Calculus of conic representations  of monotone vector fields}
\label{sec:ccr}
$\cK$-representations of pairs $(F,\X)$ ($\X$ is a nonempty convex subset of Euclidean space $E$, $F:\X\to E$ is a continuous monotone  vector field)   admit a  calculus; for verification of claims to follow, see Appendix \ref{appendix}.
\subsubsection{Raw materials}\label{raw}
Raw materials for the calculus of $\cK$-representable monotone vector fields include:
\begin{enumerate}
\item{[Affine monotone vector field]} Let $\cK$ contain all Lorentz cones. Then affine monotone vector field $F(x)=Ax+a$ on $\X=E=\bR^n$
is $\cK$-represented by conic constraint (cf. \cite[Proposition 7.4.3]{NNBook})
\begin{equation}\label{jan50}
t\geq x^T\bar{A}x+a^Tx,g=Ax+a
\end{equation}
in variables $t,g,x$, where $\bar{A}={1\over 2}[A+A^T]$.

\item{[Gradient field of continuously differentiable $\cK$-representable convex function]}  Let $\X$ be nonempty convex compact set, and $f(x):\X\to\bR$ be continuously differentiable convex function with $\cK$-representable epigraph,
\begin{equation}\label{jan61}
\{(x,s):s\geq f(x),x\in \X\}=\{(x,s):\exists u: Ax+sp+Qu\leq_{\bK} a\},
\end{equation}
the representation being essentially strictly feasible. Then the conic  constraint
\begin{equation}\label{jan60}
t\geq s+r,Ax+sp+Qu\leq_{\bK}a,r\geq a^T\lambda,A^T\lambda=g,p^T\lambda=-1,Q^T\lambda=0,\lambda\geq_{\bK^*}0
\end{equation}
in variables $t,g,x,s,r,\lambda,$ and $u$
represents $(F(\cdot):=f'(\cdot),\X)$ (cf. \cite[Proposition 7.4.4]{NNBook}).
\item{[Monotone vector field associated with continuously differentiable $\cK$-representable convex-concave function $\psi$]}  Let {$\K$ contain the 3D Lorentz cone, and let} $\U\subset\bR^{n_u}$ and $\V\subset\bR^{n_v}$ be $\cK$-representable nonempty compact sets:
\begin{equation}\label{jan70}
\begin{array}{rclr}
\U&=&\{u:\exists \alpha:\, Au+B\alpha\leq_{\bK} a\},&\qquad(a)\\
\V&=&\{v:\exists \beta: \,Cv+D\beta\leq_{\bL} b\},&(b)\\
\end{array}
\end{equation}
both representations being essentially strictly feasible. {Assume that} $\psi(u,v):\U\times \V\mapsto\bR$ {is} a continuously differentiable convex-concave function which is $\cK$-representable on $\U\times \V$ with essentially strictly feasible representation (see Section \ref{convconcreprdef}); that is, $\psi$ admits representation
\begin{equation}\label{jan71}
\forall (u\in \U,v\in \V): \psi(u,v)=\inf\limits_{f,\tau,\xi}\left\{f^Tv+\tau: Pf+\tau p+Q\xi+Ru\leq_{\bM}c\right\}
\end{equation}
such that the conic constraint $Pf+\tau p+Q\xi\leq_{\bM}c -Ru$ in variables $f,\tau,\xi$ is essentially strictly feasible for every $u\in \U$.
Let
$$
F(u,v)=[\nabla_u\psi(u,v);-\nabla_v\psi(u,v)]: \X:=\U\times\V \to E:=\bR^{n_u}\times\bR^{n_v}
$$
be the vector field associated with $\U,\V,\psi$; it is well known that this field is monotone. Then
\par\noindent
$(i)$ The set
\begin{equation}\label{jan72}
\Z=\left\{(t\in\bR,g=[h;e]\in E,x=[u;v]\in \X):\exists r,s: \begin{array}{ll}r> \max\limits_{\zeta\in\V}\left[\zeta^Te+\psi(u,\zeta)\right]&(a)\\
s>\max\limits_{\omega\in\U}\left[\omega^Th-\psi(\omega,v)\right]&(b)\\
t\geq r+s&(c)\\
\end{array}\right\}
\end{equation}
satisfies the relations
\begin{equation}\label{jan110}
\begin{array}{rclr}
\Z&\subset&\cF[F,\X]&\qquad(a)\\
x\in\X,g=F(x),t>x^TF(x)&\Rightarrow&[t;g;x]\in\Z.&(b)\\
\end{array}
\end{equation}
\par
$(ii)$ Besides this, $\cZ$
is nothing but the projection on the plane of $(t,g=[h;e],x=[u;v])$-variables of the solution set of the conic constraint
\begin{equation}\label{jan73}
\left[\begin{array}{l}
t\geq r+s,\;r\geq \theta+\gamma^Tb{+\tau},\;s\geq \theta'+c^T\delta+a^T\epsilon,\\
Au+B\alpha\leq_{\bK}a,\;Cv+D\beta\leq_{\bL} b,\;Pf+\tau p +Q\xi+Ru\leq_{\bM} c\\
C^T\gamma=f+e,D^T\gamma=0,\;P^T\delta+v=0,\;p^T\delta=-1,\;Q^T\delta=0,\;R^T\delta+A^T\epsilon=h,\;B^T\epsilon=0,\\
\gamma\geq_{\bL^*}0,\;\delta\geq_{\bM^*}0,\;\epsilon\geq_{\bK^*}0,\;
\eta\geq0,\theta\geq0,\;\eta\theta\geq1,\;
\eta'\geq0,\;\theta'\geq0,\;\eta'\theta'\geq1
\end{array}\right]
\end{equation}
in variables $t,g=[h;e],x=[u;v],r,s,f,\tau,\xi,\alpha,\beta,\gamma,\eta,\theta,\delta,\epsilon,\eta',$ and $\theta'$.\footnote{When one of the authors was a freshman at Math. Department of Moscow State University, his mate left lecture on Linear Algebra writing in his notebook: ``The lecture was terminated due to shortage of indices.'' The list of variables in (\ref{jan73}) comes close to this natural limit...}
\par\noindent
Taken together, $(i)$ and $(ii)$  say that the conic constraint  (\ref{jan73}) almost represents $(F,\X)$.
\item{[Univariate monotone rational vector field]} Let $\X=[a,b]\subset \bR$ ($-\infty<a<b<\infty$), and let $\alpha(x)$ and $\beta(x)$ be real algebraic polynomials such that $\beta(x)>0$ on $\X$. Suppose that the univariate vector field given on $\X$ by the function
$F(x)={\alpha(x)\over\beta(x)}$ is nondecreasing on $\X$, and that $\cK$ contain all semidefinite cones. Then the set
$$
\Y_+:=\Conv\{[xF(x);F(x);x], x\in\X\}+\{[t;0;0]:t\geq0\}
$$
admits explicit $\cK$-representation which, by Remark \ref{rem1}, represents $(F,\X)$.
\end{enumerate}
\subsubsection{Calculus rules}
Calculus rules are as follows.
\begin{enumerate}
\item {[Restriction on a $\cK$-representable set]} Let conic constraint
\begin{equation}\label{jan1}
Xx+Gg+tT+Uu\leq_{\bK} a
\end{equation}
in variables $t,g,x,u$
represent (almost represent) $(F,\X)$, and let set $\Y\subset E$ be $\cK$-representable:
\begin{equation}\label{jan2}
\Y=\{y\in E:\exists v: Ay+Bv\leq_{\bL} b\}
\end{equation}
and have a nonempty intersection $\Z=\Y\cap\X$ with $\X$. Denoting by $\bar{F}$ the restriction of $F$ on $\Z$, the conic constraint
\begin{equation}\label{jan3}
\{Xx+Gg+tT+Uu\leq_{\bK} a,\, Ax+Bv\leq_{\bL} b\}
\end{equation}
in variables $t,g,x,u,v$ represents (resp., almost represent) $(\bar{F},\Z)$  (cf. \cite[Proposition 7.4.5]{NNBook}).
\item {[Direct summation]} For $k\leq K$, let $\X_k$ be nonempty convex subsets of Euclidean spaces $E_k$ and $F_k:\X_k\to E_k$ be continuous monotone vector fields, and let
\begin{equation}\label{jan65}
X_kx+G_kg_k+t_k T_k+U_ku_k\leq_{\bK_k} a_k,\;1\leq k\leq K
\end{equation}
be $\cK$-representations of $(F_k,\X_k)$. Denote
\[F([x_1;...;x_K])=[F_1(x_1);...;F_K(x_K)]\;\;\mbox{and}\;\;\X=\X_1\times...\times \X_K.\]
Then the conic constraint
\begin{equation}\label{jan66}
\begin{array}{rcll}
X_kx_k+G_kg_k+t_kT_k+U_ku_k&\leq_{\bK_k}& a_k,k\leq K &(a)\\
t&=&\sum_kt_k&(b)\\
\end{array}
\end{equation}
in variables $t$, $g:=[g_1;...;g_K]$, $x:=[x_1;...;x_K]$, $t_1,...,t_K$, $u_1,...,u_K$ $\cK$-represents  $(F,\X)$.
\par
When conic constraints (\ref{jan65}) almost represent $(F_k,\X_k)$, (\ref{jan66}) almost  represents $(F,\X)$   (cf. \cite[Proposition 7.4.6]{NNBook}).

\item {[Taking conic combinations]} Let $\X$ be a nonempty convex  subset of Euclidean space $E$, {and let} $F_1,...,F_K$ be continuous monotone vector fields on $X$ with $(F_k,\X)$ admitting $\cK$-representations
\begin{equation}\label{jan4}
X_kx+G_kg+tT_k+U_ku_k\leq_{\bK_k} a_k,\;1\leq k\leq K.
\end{equation}
Let, further, $\alpha_k>0$ be given, and let
$$
F(x)=\sum_k\alpha_k F_k(x):\X\to E.
$$
The conic constraint
\begin{equation}\label{jan6}
\begin{array}{lr} X_kx+G_kg_k+t_kT_k+U_ku_k\leq_{\bK_k} a_k,\;1\leq k\leq K&\qquad(a)\\
g=\sum_k\alpha_kg_k&(b)\\
t=\sum_k\alpha_kt_k&(c)\\
\end{array}
\end{equation}
in variables $t,g,x,\{t_k,g_k,u_k,k\leq K\}$ is a $\cK$-representation of $(F,\cX)$.
\par
When (\ref{jan4}) are almost representations of $(F_k,\X)$, (\ref{jan6}) is almost representation of $(F,\X)$.

\item{[Affine substitution of variables]} Let $\X$ be a nonempty convex subset of $E$, $F$ be a continuous monotone vector field on $X$, with $(F,\X)$ given by $\cK$-representation
\begin{equation}\label{jan7}
Xx+Gg+tT+Uu\leq_{\bK} a.
\end{equation}
{Assume that}
$
\xi\mapsto A\xi+a
$
{is an} affine mapping from Euclidean space $\Lambda$ to $E$, and
$$
\Xi=\{\xi:A\xi+a\in \X\}.
$$
Then vector field $\Phi(\xi)=A^TF(A\xi+a)$: $\Xi\to\Lambda$ is continuous and monotone, and
the conic constraint
\begin{equation}\label{jan8}
\tau=t-\langle g,a\rangle ,\,\gamma=A^Tg,\;X(A\xi+a)+Gg+tT+Uu\leq_{\bK} a
\end{equation}
in variables $\tau,\gamma,\xi,g,t,$ and $u$
represents $(\Phi,\Xi)$.\\
When (\ref{jan7}) nearly represents $(F,\X)$, conic constraint (\ref{jan8}) nearly represents $(\Phi,\Xi)$.
\end{enumerate}
\subsection{Illustrations}
\subsubsection{``Academic'' illustration}
Let $\cK$ contain Lorentz cones, $M\in\bR^{n\times n}$ be such that $\overline{M}={1\over 2}[M+M^T]\succeq0$, and let $\X\subset\bR^n_+$ be nonempty and $\cK$-representable:
\begin{equation}\label{eqeq1}
\X=\{x:\exists u_0:A_0 x+B_0 u_0\leq_{\bK_0} a_0\}.
\end{equation}
{Suppose that} the operator
\begin{equation}\label{qeqreq}
F(x)=Mx-[f_1(x);...;f_n(x)]
\end{equation}
{is} monotone on $\X$,  and {that} functions $f_i(x)$ and $-\sum_ix_if_i(x)$ {are} $\cK$-representable  on $\X$:
\begin{equation}\label{qeq2}
\begin{array}{lr}
\{(x,s): x\in \X, s\geq f_i(x)\}=\{(x,s):\exists u_{1i}: A_{1i}x+sB_{1i}+C_{1i} u_{1i}\leq_{\bK_{1i}} a_{1i}\}&\qquad(a_i)\\
\{(x,s): x\in \X, s\geq -\sum_ix_if_i(x)\}=\{(x,s):\exists u_{2}: A_{2}x+sB_{2}+C_{2}u_{2}\leq_{\bK_{2}} a_{2}\}&(b)\\
\end{array}
\end{equation}
Observe, that when $M\succeq\alpha I_m$ with $\alpha>0$, $F$ definitely is monotone provided that
$f(x):=[f_1(x);...;f_n(x)]$ is Lipschitz continuous on $\X$ with Lipschitz constant w.r.t. $\|\cdot\|_2$ bounded with
$\alpha$.
\par
Given $x\in \X$ and $g\in\bR^n$, let us consider the function
$$
f_{x,g}(y)=g^Ty+F^T(y)(x-y)=-y^T{\overline M}y+\sum_i\left[g_iy_i+x_i[[My]_i-f_i(y)]+y_if_i(y)\right]:\X\to\bR.
$$
We are in the situation where $x_i\geq0$ for $x\in\X$, so that
$$
\begin{array}{ll}
&\{(y,r):y\in \X, f_{x,g}(y)\geq r\}\\
=&\left\{(y,r):\exists u_0,{s_0},s_{1i},s_{2}:\,
\begin{array}{l}
A_0y+B_0u_0\leq_{\bK_0} a_0,\;y^T\overline{M}y-s_0\leq0\\
s_{1i}+f_i(y)\leq0,\;1\leq i\leq n,s_{2}-\sum_iy_if_i(y)\leq 0\\
\sum_i[g_i+[M^Tx]_i]y_i-{s_0}+\sum_ix_is_{1i}+s_2\geq r
\end{array}
\right\}\\
=&\bigg\{(y,r):\exists u_0,{s_0},s_{1i},s_{2},u_{1i},u_{2}:\\
&\multicolumn{1}{r}{
\qquad\qquad\left.\begin{array}{l}
A_0y+B_0u_0\leq_{\bK_0} a_0,\,y^T\overline{M}y-s_0\leq0,\;A_{2}y-s_{2}B_{2}+C_{2}u_{2}\leq_{\bK_{2}}a_{2}\\
 A_{1i}y-s_{1i}B_{1i}+C_{1i}u_{1i}\leq_{\bK_{1i}} a_{1i},\;1\leq i\leq n\\
\sum_i[g_i+[M^Tx]_i]y_i{-s_0}+\sum_ix_is_{1i}+s_{2}\geq r
\end{array}
\right\}}\\
\end{array}
$$ Taking into account that $\overline{M}\succeq0$, the constraint $y^T\overline{M}y-s_0\leq0$ can be represented by strictly feasible conic constraint
$$
A_3y+s_0B_3\leq_{\bK_3} a_3
$$
on Lorentz cone ${\bL^3}$, we get
\begin{eqnarray}\label{qeq33}
\lefteqn{\max\limits_{y\in\X} f_{x,g}(y)=
\sup\limits_{y,u_0,{s_0,}s_{1i},s_{2},u_{1i},u_{2}}\bigg\{\sum_i[g_i+[M^Tx]_i]y_i{-s_0}+\sum_ix_is_{1i}+s_{2}:\nonumber}\\
&&\left.\begin{array}{l}
A_0y+B_0u_0\leq_{\bK_0} a_0,\,A_3y+s_0B_3\leq_{\bK_3} a_3\\
 A_{1i}y-s_{1i}B_{1i}+C_{1i}u_{1i}\leq_{\bK_{1i}} a_{1i},1\leq i\leq n,\,A_{2}y-s_{2}B_{2}+C_{2}u_{2}\leq_{\bK_{2i}}a_{2}\\
\end{array}\right\}.
\end{eqnarray}
Assume that the system of conic constraints (\ref{qeq33}) in variables $y,{s_0,}u_0,s_{1i},s_{2},u_{1i}$ and $u_{2}$ is essentially strictly feasible. Then, by conic quality,
\begin{eqnarray*}
\max\limits_{y\in\X} f_{x,g}(y)&=&\min_{\mu,\nu,\xi_i,\eta}\bigg\{
a_0^T\mu+a_3^T\nu+\sum_ia_{1i}^T\xi_i+a_{2}^T\eta:\nonumber\\
&&\qquad\left.\begin{array}{l}
M^Tx+g=A_0^T\mu+A_3^T\nu+\sum_iA_{1i}^T\xi_i+A_{2}^T\eta\\
x_i+B_{1i}^T\xi_i=0,\,1\leq i\leq n,\, B_{2}^T\eta=-1{,\,B_3^T\nu=-1}\\
B_0^T\mu=0,C_{1i}^T\xi_i=0,\,1\leq i\leq n,\,C_{2}^T\eta=0,\\
\mu\geq_{\bK_0^*}0,\,\nu\geq_{\bK_3^*}0,\,\xi_i\geq_{\bK_{1i}^*}0,\,1\leq i\leq n,\,\eta\geq_{\bK_{2}^*}0\\
\end{array}\right\}.
\end{eqnarray*}
Now, recalling what $f_{x,g}(y)$ is, we end up with $\cK$-representation of $\cF[F,\X]$:
\[
\begin{array}{rcl}
\cF[F,\X]&:=&\left\{[t;g;x]:x\in \X, t-g^Ty\geq F^T(y)[x-y]\,\,\forall y\in\X\right\}\\
&=&\bigg\{[t;g;x]:\exists u_0,\mu,\nu,\xi_i,\eta:\\
\multicolumn{3}{r}{\qquad\qquad\qquad\qquad\qquad\left.\begin{array}{l}
t\geq a_0^T\mu+a_3^T\nu+\sum_ia_{1i}^T\xi_i+a_{2}^T\eta,\,A_0x+B_0u_0\leq_{\bK_0}a_0\\
M^Tx+g=A_0^T\mu+A_3^T\nu+\sum_iA_{1i}^T\xi_i+A_{2}^T\eta\\
x_i+B_{1i}^T\xi_i=0,1\leq i\leq n,\, B_{2}^T\eta=-1{,\,B_3^T\nu=-1}\\
B_0^T\mu=0,C_{1i}^T\xi_i=0,1\leq i\leq n,\,C_{2}^T\eta=0\\
\mu\geq_{\bK_0^*}0,\nu\geq_{\bK_3^*}0,\xi_i\geq_{\bK_{1i}^*}0,1\leq i\leq n,\,\eta\geq_{\bK_{2}^*}0
\end{array}\right\},}
\end{array}
\]
so that the conic constraint
\[
\left[\begin{array}{l}
t\geq a_0^T\mu+a_3^T\nu+\sum_ia_{1i}^T\xi_i+a_{2}^T\eta,\,
A_0x+B_0u_0\leq_{\bK_0}a_0\\
x_i+B_{1i}^T\xi_i=0,1\leq i\leq n,\, B_{2}^T\eta=-1,{,\,B_3^T\nu=-1}\\
B_0^T\mu=0,C_{1i}^T\xi_i=0,1\leq i\leq n,\,C_{2}^T\eta=0\\
\mu\geq_{\bK_0^*}0,\nu\geq_{\bK_3^*}0,\xi_i\geq_{\bK_{1i}^*}0,1\leq i\leq n,\,\eta\geq_{\bK_{2}^*}0
\end{array}\right]
\]
in variables $t,g,x,u_0,\mu,\nu,\xi_i$ and $\eta$
is a $\cK$-representation of $(F,\X)$.
\subsubsection{Nash Equilibrium}\label{Ill2}
The ``covering story'' for {this example} is as follows:
\begin{quote}
$n\geq2$ retailers intend to enter certain market, say, one of red herrings. To this end they should select their ``selling capacities'' (say, rent areas at malls) $x_i$, $1\leq i\leq n$, in given ranges $\cX_i=[0,X_i]$ ($0<X_i<\infty$). With the selections $x=[x_1;...;x_n]\in\X=\X_1\times...,\times \X_n\subset\bR^n_+$, the monthly loss of the $i$-th retailer is
$$
\phi_i(x)=c_ix_i-{x_i\over\sum_{j=1}^Kx_j+a}b,
$$
where $c_ix_i$, $c_i>0$, is the price of the capacity $x_i$, $b>0$ is the dollar {value of the} demand on red herrings, and $a>0$ is the total selling  capacity of the already existing retailers; the term $-{x_i\over\sum_{j=1}^nx_j+a}b$ is the minus {the} revenue of the $i$-th retailer under the assumption that the total demand is split between retailers proportionally to their selling capacities.
\par
We want to find is a Nash equilibrium -- a point $x^*=[x^*_1;...;x^*_n]\in\X$ such that every one of the functions $\phi_i(x_i^*,...,x_{i-1}^*,x_i,x_{i+1}^*,...,x_n^*)$ attains its minimum over $x_i\in\X_i$ at $x_i=x_i^*$, so that for the $i$-th retailer there is no incentive to deviate from capacity $x_i^*$ provided that all other retailers $j$ stick to capacities $x_j^*$, and this is so for all $i$.
\end{quote}
As is immediately seen, the Nash Equilibrium problem in question is convex, meaning that $\phi_i(x)$ are convex in $x_i$ and concave in $x^i=(x_i,...,x_{i-1},x_{i+1},...,x_n)$ and, on the top of it, $\sum_i\phi_i(x)$ is convex on $\X$. It is well known (for justification, see, e.g., \cite{NOR}) that for such a problem, the vector field
$$
F(x)=\left[{\partial \phi_1(x)\over\partial x_1};{\partial \phi_2(x)\over\partial x_2};...;{\partial \phi_n(x)\over\partial x_n}\right]:\X\to\bR^n
$$ is monotone, and that Nash Equilibria are exactly the weak$=$strong solutions to $\VI(F,\X)$. Specifying $\cK$ as the family of direct products of Lorentz cones, we are about to demonstrate that $F$ admits an explicit $\cK$-representation on $\X$, allowing to reduce the problem of finding Nash Equilibrium to an explicit Second Order conic problem. This is how the construction goes (to save notation, in what follows we set ${b}=a=1$, which always can be achieved by appropriate scaling of the capacities and loss functions):
\begin{enumerate}
\item Observe that all we need is a $\cK$-representation {$\R(t,g,x,u)$} of $(F,\X)$---a $\cK$-conic constraint in variables $t\in\bR$, $g\in\bR^n$, $x\in \bR^n$ and additional variables $u$---satisfying the requirements specified in Section \ref{cKrMVF}. By Proposition \ref{propreprjan}, given $\R$, for $\epsilon\geq0$ we can write down an explicit $\cK$-conic constraint $\cC_\epsilon(t,g,x,u,v)$ in variables $t,g,x$ and additional variables $u,v$,
    with the size of the constraint (dimension of the associated cone and the total number of variables) independent of $\epsilon$ such that the constraint is feasible, and the $x$-component of any  feasible solution to the constraint is an $\epsilon$-solution to $\VI(F,\X)$:
$$
x\in \X\ \&\ \epsVI[x|F]=\max_{y\in \X} F^T(y)[x-y]\leq\epsilon.
$$
As a result, finding $\epsilon$-solution to $\VI(F,\X)$ reduces to solving an explicit solvable Second Order feasibility conic problem of the size independent of $\epsilon$.
\item Thus, all we need is to find an explicit $\cK$-representation of $(F,\X)$. This can be done as follows:
\begin{enumerate}
\item Consider the convex-concave function
$$
\psi(u,v)=-{u\over u+v+1}:\bR^2_+\to\bR
$$
along with the associated monotone vector field
$$
\Phi(u,v)=\left[{\partial \psi(u,v)\over\partial u};-{\partial \psi(u,v)\over\partial v}\right]=\left[-{v+1\over(u+v+1)^2};-{u\over(u+v+1)^2}\right].
$$
Let also $A_i$ be $2\times n$ matrix with $i$-th column $[1;0]$ and all remaining columns equal to $[0;1]$, and let $G(x)=\nabla\left({1\over \sum_jx_j+1}\right)$. As is immediately seen, we have
$$
\forall x\in\bR^n_+:\Psi(x):= 2c+\sum_{i=1}^nA_i^T\Phi(A_ix)+G(x)=2F(x),
$$
so that a $\cK$-representation of $F$ is readily given by  $\cK$-representations of the constant monotone vector field $\equiv c$, of $\Phi$, and of $G$ via our calculus (rules on affine substitution of
argument and on summation). \par
A $\cK$-representation of the constant vector field  $\equiv c$ is trivial; it is given by the system of linear equalities $t=c^Tx,g=c$ on variables $(t,g,x)\in\bR\times\bR^n\times\bR^n$.
A $\cK$-representation of the gradient vector field $G(x)$ is given by calculus rule 2 in Section \ref{raw}; we are in the case where (\ref{jan61}) reads
$$
\begin{array}{rcl}
\{(x,s): s\geq f(x),x\in X\}&=&\bigg\{(x,s): -x\leq 0,x\leq X:=[X_1;...;X_n],\\
&&\quad\underbrace{-\big[0;\hbox{\footnotesize$\sum$}_ix_i-s;\hbox{\footnotesize$\sum$}_ix_i+s\big]\leq_{\bL^3}[2;1;1]}_{
\hbox{\footnotesize$\begin{array}{ll}\Leftrightarrow &x+1+s\geq0,\\
&(\hbox{\footnotesize$\sum$}_ix_i+1+s)^2\geq (\hbox{\footnotesize$\sum$}_ix_i+1-s)^2+4\\
\Leftrightarrow& s(\sum_ix_i+1)\geq1\\
\end{array}$}}\bigg\},
\end{array}
$$
so that (\ref{jan60}), as seen from immediate computation, reduces to the system of constraints
$$
\begin{array}{c}
t\geq r+s,\, 0\leq x\leq X,\,r\geq \sum_i\max[X_i(g_i+\theta),0]+\theta-2\nu\\
\underbrace{\big[2;\hbox{\footnotesize$\sum$}_ix_i+1-s;\hbox{\footnotesize$\sum$}_ix_i+1+s\big]\geq_{\bL^3}0}_{\Leftrightarrow s(\sum_jx_j+1)\geq1},\,\underbrace{[2\nu;\theta-1;\theta+1]\geq_{\bL^3}0}_{\Leftrightarrow \nu^2\leq\theta}\\
\end{array}$$
in variables  $t,g,x,\theta,\nu$.

\item It remains to find a $\cK$-representation of the vector field $\Phi$ on a rectangle $\Xi=\U\times\V$, $\U=[0,U]$, $\V=[0,V]$ with given $U>0$, $V>0$ (they should be large enough in order for $A_ix$, $x\in\X$, to take values in the rectangle).

    \par
    Let us use the construction described in item {3} of Section \ref{raw}.  By (i) of item 3, the set
$$
\overline{\Z}=\left\{(t\in\bR,g=[h;e]\in \bR^2,x=[u;v]\in \Xi):\exists r,s: \begin{array}{ll}r\geq \max\limits_{\zeta\in\V}\left[\zeta e+\psi(u,\zeta)\right]\\
s\geq\max\limits_{\omega\in\U}\left[\omega^Th-\psi(\omega,v)\right]\\
t\geq r+s\\
\end{array}\right\}
$$
satisfies the relations
$$
\begin{array}{rcl}
\overline{\Z}&\subset&\cF[\Phi,\Xi]\\
x\in\X,g=\Phi(x),t\geq x^T\Phi(x)&\Rightarrow&[t;g;x]\in\overline{\Z},
\end{array}
$$
so that  all we need to get a $\cK$-representation of $\Phi$ is to represent $\overline{\Z}$ as a projection onto the plane of $t,g,x$-variables of the solution set of $\cK$-conic constraint in variables $t,g,x$ and additional variables. To this end note
that for all $(u,v)\in\Xi=\U\times \V$ one has
    \newcommand{\bBigl}[1]{\left\{\rule{0 cm}{#1}\right.}
\newcommand{\bBigr}[1]{\left.\rule{0 cm}{#1}\right\}}
\begin{equation}\label{observe1}
{
\psi(u,v)=
\min\limits_{f,t,s,\tau}\bBigl{1cm} fv+t:\,\underbrace{\left[\begin{array}{l} \tau\geq0,
\,1\geq f\geq0,\,s\geq0,\,t\leq 1,\,s+\tau\leq 1\\
{[2s;u-f;u+f]\in\bL^3}\\
{[2(1-s);t-f;t-f+2]\in\bL^3}\\
{[2;u-\tau+1;u+\tau+1]\in\bL^3}
\end{array}\right]}_{\Leftrightarrow \,(f,t,s,\tau)\in\Pi_+}\bBigr{1cm}
}
\end{equation}
and
\begin{equation}\label{observe2}
{
-\psi(u,v)
=\min\limits_{f,t,s}\bBigl{0.8cm}fu+t:
\,\underbrace{\left[\begin{array}{l} 1\geq f\geq0,\,0\leq s\leq 1,\,t\leq 1\\
{[2s;v+1-f;v+1+f]\in\bL^3}\\
{[2(1-s);t-1;t+1]\in\bL^3}
\end{array}\right]}_{\Leftrightarrow \,(f,t,s)\in\Pi_-}\bBigr{0.8cm},
}
\end{equation}
see Section \ref{justif} for justification.
{Since $\Pi_\pm$ are nonempty convex compact sets, by the Sion-Kakutani Theorem one has from (\ref{observe1})}:
\[
\begin{array}{rcl}
\max\limits_{\zeta\in\V}\left[\zeta e+\psi(u,\zeta)\right]&=&\max\limits_{0\leq\zeta\leq V}\min\limits_{f,t,s,\tau}\left\{[f+e]\zeta+t:(f,t,s{,\tau})\in\Pi_+\right\}\\
&=&\min\limits_{f,t,s{,\tau}}\left\{\max\limits_{0\leq\zeta\leq V}[f+e]\zeta+t:(f,t,s{,\tau})\in\Pi_+\right\}\\
&=&\min\limits_{f,t,s{,\tau}}\left\{\max[0,V(f+e)]+t:(f,t,s{,\tau})\in\Pi_+\right\},
\end{array}
\] {and from (\ref{observe2}):}
\[
\begin{array}{rcl}\max\limits_{\omega\in\U}\left[\omega h-\psi(\omega,v)\right]&=&\max\limits_{0\leq\omega\leq U}\min\limits_{f,t,s}\left\{[f+h]\omega+t:(f,t,s)\in\Pi_-\right\}\\
&=&\min\limits_{f,t,s}\left\{\max\limits_{0\leq\omega\leq U}[f+h]\omega+t:(f,t,s)\in\Pi_-\right\}\\
&=&\min\limits_{f,t,s}\left\{\max[0,U(f+h)]+t:(f,t,s)\in\Pi_-\right\}.
\end{array}
\]
As a result,
\[
\begin{array}{l}
\overline{\Z}=\bigg\{(t\in\bR,g=[h;e]\in \bR^2,0\leq x:=[u;v]\leq [U;V]):\exists r,s,f_\pm,t_\pm,s_\pm{,\tau}:\\
\multicolumn{1}{r}{\begin{array}{l}
t\geq r+s,
r\geq \max[0,V(f_++e)]+t_+,
s\geq\max[0,U(f_-+h)]+t_-\\
\,(f_+,t_+,s_+{,\tau})\in\Pi_+,\,(f_-,t_-,s_-)\in\Pi_-\end{array}\bigg\}.}
\end{array}
\]
Recalling what are $\Pi_\pm$, {we can straightforwardly} represent $\overline \Z$ as projection onto the space of $t,g,x$-variables of a set given by a $\cK$-conic inequality, {ending} up with an explicit $\cK$-representation of $\Phi$.
\end{enumerate}
\end{enumerate}
{\begin{remark} {\rm The just outlined construction can be {used in the case of} a general convex Nash Equilibrium problem where, given $n$ convex compact sets $\X_i\subset\bR^{k_i}$ and $n$ continuously differentiable functions
$\phi_i(x_1,...,x_n): \,\X:=\X_1\times...\times \X_n\to\bR$ with $\phi_i$ convex in $x_i$, concave in $(x_1,...,x_{i-1},x_{i+1},...,x_n)$ and such that $\sum_i\phi_i$ is convex, one is looking for Nash Equilibria---points $x^*=[x^*_1;...;x^*_n]\in\X$ such that for every $i$, the function $\phi_i(x_1^*,...,x_{i-1}^*,x_i,x_{i+1}^*,...,x_n^*)$ attains its minimum over $x_i\in \X_i$ at $x_i=x_i^*$. As it was already mentioned, these equilibria are exactly the weak$=$strong solutions to $\VI(F,\X)$ where the monotone vector field $F$ is given by $F(x_1,...,x_n)=\left[{\partial\phi_1(x)\over\partial x_1};...;{\partial\phi_n(x)\over\partial x_n}\right]$.
Observing that the monotone vector fields
$$F^i(x)=\left[-{\partial\phi_i(x)\over\partial x_1};...;-{\partial\phi_{i}(x)\over\partial x_{i-1}};{\partial\phi_i(x)\over\partial x_i};-{\partial\phi_{i}(x)\over\partial x_{i+1}};...;-{\partial\phi_i(x)\over\partial x_n}\right]$$
associated with convex-concave functions $\phi_i$
and the monotone vector field $F^0(x)=\nabla\left(\sum_i\phi_i(x)\right)$ are linked to $F$ by the relation $2F=\sum_{i=0}^n F^i$, we see that in order to get a $\cK$-representation of $F$, it suffices to have at our disposal $\cK$-representations of $F^i$, $0\leq i\leq n$. These latter representations, in good cases, can be built according to recipes presented in items 2 and 3 of Section \ref{raw}.}
\end{remark}}\noindent
The latter remark puts in proper perspective the ``red herring'' illustration which by itself is of no actual interest: setting $s=\sum_ix_i+1$,  $x\in\X$ is a Nash Equilibrium if and only if
$$
c_i-{1\over s}+{x_i\over s^2}={\partial\phi_i(x)\over\partial x_i}\left\{\begin{array}{ll}\geq0,&x_i=0\\
=0,&0<x_i<X_i\cr
\leq0,&x_i=X_i.
\end{array}\right.
$$
As a result, finding the equilibrium reduces to solving on the segment $s\in\big[1,\sum_iX_i+1\big]$ the univariate equation
$$
\sum_ix_i(s)+1=s, \,x_i(s)=\left\{\begin{array}{ll}0,&s(1-sc_i)< 0\\
s(1-sc_i),&0\leq s(1-sc_i)
\leq X_i\\
X_i,&s(1-sc_i)>X_i.
\end{array}\right.
$$
\section*{Acknowledgment}
This work was supported by Multidisciplinary Institute in Artificial intelligence MIAI {@} Grenoble Alpes (ANR-19-P3IA-0003).

\appendix
\section{Derivations for Section \protect{\ref{sec:ccr}}}\label{appendix}
\subsection{Verification of ``raw materials''}
\begin{enumerate}\item{\em [Affine monotone vector field]}
Note first that the symmetric part $\bar{A}$  of $A$ is $\succeq0$ due to the monotonicity of $F$ on $E$, so that (\ref{jan50}) can be rewritten as a conic constraint on a direct product of properly selected Lorentz cone and nonnegative orthant, and this direct product belongs to $\cK$. Second, when $t,g,x$ solve (\ref{jan50}), for every $y\in E$ one has $x^TAx=x^T\bar{A}x$, whence
$$
t-y^Tg=t-y^T[Ax+a]=\underbrace{t-x^T[Ax+a]}_{=t-x^T\bar{A}x-x^Ta\geq0}-[y-x]^T[Ax+a]\geq [x-y]^T[Ax+a]\geq [x-y]^T[Ay+a],
$$
that is, $[t;g;x]\in\cF[F,E]$, as required in item (i) of a representation of $(F,E)$.  Furthermore, when $x\in E$, setting $g=F(x)=Ax+a$ and $t=x^TF(x)=x^T[Ax+a]=x^T[\bar{A}x+a]$, we see that $[t;g;x]$ satisfy (\ref{jan50}), as required in item (ii) of a representation of $(F,E)$.
\item{\em [Gradient field of continuously differentiable $\K$-representable convex function]}
Let \[f_*(y)=\sup\limits_{x\in \X}[x^Ty-f(x)]\] be the Fenchel transform of $f$:
\[
\begin{array}{rcl}
f_*(y)&=&\sup\limits_{x\in \X}[x^Ty-f(x)]=\sup\limits_{x\in\X,\;s\geq f(x)}[x^Ty-s]\\
&=&\sup\limits_{x,u}[x^Ty-s: Ax+sp+Qu\leq_{\bK} a]\\
&=&\min\limits_{\lambda}\left\{a^T\lambda: A^T\lambda=y,p^T\lambda=-1,Q^T\lambda=0,\lambda\geq_{\bK^*} 0\right\}\;\;\hbox{[by conic duality],}
\end{array}
\]
that is,
\begin{equation}
\label{jan62}
r\geq f_*(g)\;\Leftrightarrow\;\exists \lambda:r\geq a^T\lambda, A^T\lambda=g,p^T\lambda=-1,Q^T\lambda=0,\lambda\geq_{\bK^*} 0.
\end{equation}
Next, let
\begin{equation}
\label{jan162}
\begin{array}{rclr}
\Z&=&\{(t,g,x): \exists s,r: \;s\geq f(x),r\geq f_*(g),s+r\leq t\}&\qquad(a)\\
&=&\bigg\{(t,g,x):\exists s,r,u,\lambda: \;\begin{array}{l}t\geq s+r,Ax+sp+Qu\leq_{\bK}a\\
r\geq a^T\lambda,A^T\lambda=g,p^T\lambda=-1,Q^T\lambda=0,\lambda\geq_{\bK^*}0\end{array}\bigg\}&(b)\\
\end{array}
\end{equation}
where $(\ref{jan162}.b)$ is due to (\ref{jan61}) and (\ref{jan62}). Thus, $\Z$ is the projection of the solution set of (\ref{jan60}) on the space of $t,g,x$-variables, and we should check that
\par
$\quad$(i) $\cZ\subset\cF[f'(\cdot),\X]$, and
\par
$\quad$(ii) when $x\in \X$, $g=f'(x)$, and $t=x^Tf'(x)$, we have $[t;g;x]\in\cZ$.
\par\noindent
(i): Let $(t,g,x)\in\Z$. We have $t\geq s+r$ for properly selected $s\geq f(x)$ and $r\geq f_*(g)$, whence for $y\in \X$ it holds
$$
t-g^Ty\geq f(x)+f_*(g)-g^Ty\geq f(x)-f(y)\geq (x-y)^Tf'(y),
$$
that is, $[t;g;x]\in\cF[f'(\cdot),\X]$. (i) is proved.\\
(ii):  Given  $x\in \X$, let us set $g=f'(x)$, $t=x^Tg$, $s=f(x)$, $r=f_*(g)$. We have \[r=f_*(g)=\sup\limits_{x'} \{g^Tx'-f(x'):\,x'\in\X\}=g^Tx-f(x)
\] (the concluding equality is due to $g=f'(x)$); thus, $t=r+s$. Invoking  $(\ref{jan162}.a)$, we conclude that
$(t,g,x)\in \Z$. (ii) is proved.
\item{[Monotone vector field associated with $\cK$-representable convex-concave function $\psi$]}~
\vspace{-0.5cm}\paragraph{Verifying $(i)$:}~
\\{\bf 1$^o$.}
Let us show (\ref{jan110}.$a$). Let $(t,g=[h;e],x=[u;v])\in\Z$, and let $r,s$ be reals which, taken together with $t,g,x$, form a feasible solution to the system of constraints in (\ref{jan72}). For $y=[w;z]\in\X$ we have
$$
\begin{array}{rcl}
t-y^Tg&\geq&r+s-w^Th-z^Te\;\;\hbox{[by (\ref{jan72}.$c$)]}\\
&\geq&\left[z^Te+\psi(u,z)\right] +\left[w^Th-\psi(w,v)\right] -w^Th-z^Te\hbox{\ [by (\ref{jan72}.$a$,$b$)]}\\
&=&\psi(u,z)-\psi(w,v)=[\psi(u,z)-\psi(w,z)]+[\psi(w,z)-\psi(w,v)]\\
&\geq&[u-w]^T\nabla_u\psi(w,z)-[v-z]^T\nabla_v\psi(w,z)\;\;\hbox{[because $\psi$ is convex-concave]}\\
&=&[x-y]^TF(y),\\
\end{array}
$$
implying that $\Z\subset \cF[F,\X]$.
\\{\bf 2$^o$.}
Let us now verify (\ref{jan110}.$b$). Given $x=[u;v]\in \X$, let us set
$$
g=F(x)=[h;e],\, h=\nabla_u\psi(u,v),\,e=-\nabla_v\psi(u,v);\,\bar{t}=x^TF(x)=h^Tu+e^Tv,
$$
and let $t>\bar{t}$.
The function $\zeta^Te+\psi(u,\zeta)$ of $\zeta\in\V$ is concave, and its gradient w.r.t. $\zeta$ taken at the point $\zeta=v$ vanishes, implying that
$$
\bar{r}:=v^Te+\psi(u,v)\geq \zeta^Te+\psi(u,\zeta)\,\,\forall \zeta\in\V.
$$
Similarly,  the function $\omega^Th-\psi(\omega,v)$ of $\omega\in\U$ is concave. and its gradient w.r.t. $\omega$ taken at the point $\omega=u$ vanishes, implying that
$$
\bar{s}:=u^Th-\psi(u,v)\geq \omega^Th-\psi(\omega,v)\,\,\forall \omega\in\U.
$$
Observing that $\bar{r}+\bar{s}=\bar{t}<t$, we can find $r>\bar{r}$ and $s>\bar{s}$ in such a way that $r+s\leq t$. Looking at the constraints (\ref{jan72}.$a$-$c$), we conclude that whenever $x=[u;v]\in \X$ and $t>x^TF(x)$, the triple $(t,g=F(x),x)$ can be augmented by $r,s$ to satisfy all constraints (\ref{jan72}.$a$-$c$), that is, $(t,g,x)\in\Z$, as claimed in (\ref{jan110}.$b$).
\paragraph{Verifying $(ii)$:}
Clearly, all we need in order to justify claim in $(ii)$ is to show that the set
$$
\Z^+=\left\{(g=[h;e],x=[u;v],r,s): u\in\U,v\in\V, r> \max\limits_{z\in\V}\left[z^Te+\psi(u,z)\right],s>\max\limits_{w\in\U}\left[w^Th-\psi(w,v)\right] \right\}
$$
is the projection of the solution set of the conic constraint
\begin{equation}\label{jan80}
\left[\begin{array}{l}
r\geq \theta+\gamma^Tb{+\tau},\;\eta\geq0,\;\theta\geq0,\;\eta\theta\geq1\\
s\geq \theta'+ c^T\delta+a^T\epsilon,\eta'\geq0,\theta'\geq0,\eta'\theta'\geq1\\
Au+B\alpha\leq_{\bK}a,\;Cv+D\beta\leq_{\bL} b,\;Pf+\tau p +Q\xi+Ru\leq_{\bM} c\\
C^T\gamma=f+e,D^T\gamma=0,\gamma\geq_{\bL^*}0,;\delta\geq_{\bM^*}0,\;\epsilon\geq_{\bK^*}0\\
P^T\delta+v=0,p^T\delta=-1,Q^T\delta=0,R^T\delta+A^T\epsilon=h,B^T\epsilon=0
\end{array}\right]
\end{equation}
in variables $x=[u;v],g=[h;e],r,s,\alpha,\beta,\gamma,f,\tau,\xi,\delta,\epsilon,\eta,\theta,\eta',$ and $\theta'$ on the plane of $g=[h;e],x=[u;v],r,s$-variables. Here is the proof.
\\{\bf 1$^o$.} Recall that $\V$ is convex and compact. Thus, we have
$$
\begin{array}{l}
\forall (x=[u;v]\in\X,g=[h;e]):\\
\max\limits_{z\in \V}\left[z^Te+\psi(u,z)\right]
=\max\limits_{z\in\V} \left[z^Te+\inf\limits_{f,\tau,\xi}\left[f^Tz+\tau:Pf+\tau p +Q\xi+Ru\leq_{\bM} c\right]\right]  \hbox{\ [by (\ref{jan71})]}\\
=\max\limits_{z\in\V} \inf\limits_{f,\tau,\xi}\left[(f+e)^Tz+\tau:Pf+\tau p +Q\xi+Ru\leq_{\bM} c\right]\\
= \inf\limits_{f,\tau,\xi}\left\{\max\limits_{z\in\V}(f+e)^Tz+\tau:Pf+\tau p +Q\xi+Ru\leq_{\bM} c\right\} \ \left[\hbox{by the Sion-Kakutani Theorem}\right]\\
=\inf\limits_{f,\tau,\xi}\left[\max\limits_{z,\beta}\left[(f+e)^Tz+\tau: Cz+D\beta\leq_{\bL} b\right]:Pf+\tau p +Q\xi+Ru\leq_{\bM} c\right] \hbox{\ [by (\ref{jan70}.$b$)]}\\
=\inf\limits_{f,\tau,\xi,\gamma}\left\{\gamma^Tb{+\tau}: C^T\gamma=f+e,D^T\gamma=0,\gamma\geq_{\bL^*}0,Pf+\tau p +Q\xi+Ru\leq_{\bM} c\right\}\\
\hbox{[by conic duality; recall that (\ref{jan70}.$b$) is essentially strictly feasible].}\\
\end{array}
$$
Together with (\ref{jan70}) the latter relation results in
\begin{eqnarray}\label{jan49}
\lefteqn{\big\{(x=[u;v],g=[h;e],r):\,x\in\X, r>\max\limits_{z\in \V}\left[z^Te+\psi(u,z)\right]\big\}}\nonumber\\
&=\bigg\{(x=[u;v],g=[h;e],r):&\exists f,\tau,\xi,\alpha,\beta,\gamma:\nonumber\\
&&\left.\begin{array}{l}r>\gamma^Tb{+\tau},\;Au+B\alpha\leq_{\bK}a,\;Cv+D\beta\leq_{\bL} b,\;\gamma\geq_{\bL^*}0\\
C^T\gamma=f+e,\;D^T\gamma=0,\;Pf+\tau p +Q\xi+Ru\leq_{\bM} c
\end{array}\right\}\nonumber\\
&=\bigg\{(x=[u;v],g=[h;e],r):&\exists f,\tau,\xi,\alpha,\beta,\gamma,\eta,\theta:\nonumber\\
&&\left.\begin{array}{l}r\geq\gamma^Tb{+\tau}+ \theta,\;\eta\geq0,\;\theta\geq0,\;\eta\theta\geq1\\
Au+B\alpha\leq_{\bK}a,\;Cv+D\beta\leq_{\bL} b,\;,\gamma\geq_{\bL^*},\\
C^T\gamma=f+e,\;D^T\gamma=0,\;Pf+\tau p +Q\xi+Ru\leq_{\bM} c
\end{array}\right\}.
\end{eqnarray}
\\{\bf 2$^o$.} When $v\in\V$ and $h\in\bR^{n_u}$ we have
$$
\begin{array}{l}
\max\limits_{w\in\U}\left[w^Th-\psi(w,v)\right]=\max\limits_{w\in\U}\left[w^Th+\sup\limits_{f,\tau,\xi}\left[-f^Tv-\tau:Pf+\tau p +Q\xi+Rw\leq_{\bM} c\right]\right]
\hbox{\ [by (\ref{jan71})]}\\
=\sup\limits_{w,f,\tau,\xi,\alpha}\left[w^Th-f^Tv-\tau:Pf+\tau p +Q\xi+Rw\leq_{\bM} c,\;Aw+B\alpha\leq_{\bK} a\right] \hbox{\ [by (\ref{jan70}.$a$)]}\\
=\min\limits_{\delta,\epsilon}\bigg\{c^T\delta+a^T\epsilon:\begin{array}{l}P^T\delta+v=0,\;p^T\delta=-1,Q^T\delta=0,\;R^T\delta+A^T\epsilon=h\\
B^T\epsilon=0,\;\delta\geq_{\bM^*}0,\;\epsilon\geq_{\bK^*}0\end{array}\bigg\}\\
\left[\hbox{\begin{tabular}{l}by conic duality; recall that the $\cK$-representations (\ref{jan70}.$a$)\\
 of $\U$ and (\ref{jan71}) of $\psi$ are essentially strictly feasible.
 \end{tabular}}\right]
\end{array}
$$
Taken together with (\ref{jan70}), the latter relation results in
\begin{eqnarray}\label{jan48}
\lefteqn{\big\{(x=[u;v],g=[h;e],s): x\in\X, s>\max\limits_{w\in \U}\left[h^Tw-\psi(w,u)\right]\big\}}\nonumber\\
&=&\bigg\{(x=[u;v],g=[h;e],s):\exists \alpha,\beta,\delta,\epsilon:\nonumber\\
&&\qquad\qquad\left.\begin{array}{l}s>c^T\delta+a^T\epsilon,\;Au+B\alpha\leq_{\bK}a,\;Cv+D\beta\leq_{\bL} b,\;\delta\geq_{\bM^*}0,\;\epsilon\geq_{\bK^*}0\\
P^T\delta+v=0,\;p^T\delta=-1,\;Q^T\delta=0,\;R^T\delta+A^T\epsilon=h,\;B^T\epsilon=0
\end{array}\right\}\nonumber\\
&=&\bigg\{(x=[u;v],g=[h;e],s):\exists \alpha,\beta,\delta,\epsilon,\eta',\theta':\nonumber\\
&&\qquad\qquad\left.\begin{array}{l}s\geq \theta'+c^T\delta+a^T\epsilon,\;\eta'\geq0,\;\theta'\geq0,\;\eta'\theta'\geq1\\
Au+B\alpha\leq_{\bK}a,\;Cv+D\beta\leq_{\bL} b,\;\delta\geq_{\bM^*}0,\;\epsilon\geq_{\bK^*}0\\
P^T\delta+v=0,\;p^T\delta=-1,\;Q^T\delta=0,\;R^T\delta+A^T\epsilon=h,\;B^T\epsilon=0
\end{array}\right\}.
\end{eqnarray}
Finally, (\ref{jan49}) and (\ref{jan48}) together imply that
$$
\begin{array}{l}
\Z^+:=\bigg\{(g=[h;e],x=[u;v],r,s): u\in\U,v\in\V,
{ r> \max\limits_{z\in\V}\left[z^Te+\psi(u,z)\right],s>\max\limits_{w\in\U}\left[w^Th-\psi(v,u)\right]\bigg\}}\\
=\bigg\{(g=[h;e],x=[u;v],r,s):\exists f,\tau,\xi,\alpha,\beta,\gamma,\eta,\theta,\delta,\epsilon,\eta',\theta':\\
\multicolumn{1}{r}{\qquad\qquad\left.
\begin{array}{l}
r\geq \theta+\gamma^Tb{+\tau},\eta\geq0,\theta\geq0,\eta\theta\geq1\\
s\geq \theta'+c^T\delta+a^T\epsilon,\eta'\geq0,\theta'\geq0,\eta'\theta'\geq1\\
Au+B\alpha\leq_{\bK}a,\;Cv+D\beta\leq_{\bL} b,\;Pf+\tau p +Q\xi+Ru\leq_{\bM} c,\gamma\geq_{\bL^*}0,\delta\geq_{\bM^*}0,\epsilon\geq_{\bK^*}0\\
C^T\gamma=f+e,D^T\gamma=0,
P^T\delta+v=0,p^T\delta=-1,Q^T\delta=0,R^T\delta+A^T\epsilon=h,B^T\epsilon=0\\
\end{array}\right\},}
\end{array}
$$
as claimed in (\ref{jan80}).
\item{[Univariate monotone rational vector field]} For evident reasons, it suffices to consider the case of $\X=[0,1]$; recall that $\beta(t)>0$ on $\X$. Let degrees of $\alpha$ and $\beta$ be $\mu$ and $\nu$, respectively, and let
$$
\kappa=\max[\mu,\nu]+1.
$$
\paragraph{1$^0$.}
Consider the curves
$$
\delta(t)=[tF(t);F(t);t]={1\over \beta(t)} [t\alpha(t);\alpha(t);t\beta(t)]:[0,1]\to\bR^3,\,\gamma(t)={1\over\beta(t)}\left[1;t;t^2;...;t^\kappa\right]:[0,1]\to\bR^{\kappa+1}
$$
For properly selected matrix $A$ we have
$$
\delta(t)=A\gamma(t),\,0\leq t\leq 1,
$$
whence
$$
\Y:=\Conv\{\delta(t),0\leq t\leq 1\}=A\Z,\,\Z=\Conv\{\gamma(t),0\leq t\leq 1\}.
$$
We intend to build a semidefinite representation (SDR) of $\Y$ (i.e., $\cK$-representation with $\cK$ comprised of finite direct products of semidefinite cones). Semidefinite representability (and $\cK$-representability in general) of a set is preserved when taking linear images: an SDR
$$
\Z=\{z:\exists u: \cA(z,u)\succeq0\}\eqno{[\hbox{$\cA(z,u)$ is affine in $[z;u]$ symmetric matrix]}}
$$
of $\Z$ implies the representation
$$
\Y:=A\Z=\{y:\exists [z;u]: y=Az,\cA(z,u)\succeq0\},
$$
and the system of the right hand side constraints can be written down as a Linear Matrix Inequality in variable $y$ and additional variables $z,u$. Thus, all we need is to build an SDR of $\Z$.
\paragraph{2$^0$.} We shall get SDR of $\cZ$ from SDR of the ``support cone''
$$
\P=\{[p;q]\in\bR^{\kappa+1}\times\bR: \min_{z\in\Z} p^Tz-q\geq0\}=\{[p;q]:\min_{0\leq t\leq 1}p^T\gamma(t)\geq q\}
$$
of $\cZ$.\par
Given $p=[p_0;p_1;...;p_{\kappa+1}]\in\bR^{\kappa+1}$, let, with a slight abuse of notation, $p(t)=\sum_{i=0}^\kappa p_it^i$ be the polynomial with coefficients $p_i$, $0\leq i\leq\kappa$.
We have
$$
\begin{array}{l}
[p;q]\in\cP\Leftrightarrow {p(t)\over \beta(t)}\geq q\,\forall t\in[0,1]\\\Leftrightarrow {(1+\tau^2)^\kappa p(\tau^2/(1+\tau^2))\over (1+\tau^2)^\kappa\beta(\tau^2/(1+\tau^2))}\geq q\;\forall \tau\in\bR\\
\Leftrightarrow \underbrace{(1+\tau^2)^\kappa p(\tau^2/(1+\tau^2))-q(1+\tau^2)^\kappa\beta(\tau^2/(1+\tau^2))}_{=:\pi_{p,q}(\tau)}\geq0\;\forall \tau\in\bR.
\end{array}
$$
Note that $\pi_{p,q}(\tau)$ is a polynomial of $\tau$ of degree $\leq2\kappa$, and the vector $\pi_{p,q}$ of coefficients of polynomial is linear in $[p;q]$: $\pi_{p,q}=P[p;q]$. We see that $\cP$ is the inverse image
of the cone $\cP_{2\kappa}$ of coefficients of polynomials of degree $\leq 2\kappa$ which are nonnegative on the entire axis:
$$
\cP=\{[p;q]:P[p;q]\in\cP_{2\kappa}\}.
$$
As was observed by Yu. Nesterov \cite{Nesterov2000}, the cone $\cP_{2\kappa}$ is the linear image of the semidefinite cone $\bS^{\kappa+1}_+$:
$$
\cP_{2\kappa}=\{\pi\in\bR^{2\kappa+1}:\exists x=[x_{ij}]_{0\leq i,j\leq\kappa}\in\bS^{\kappa+1}_+: [\cQ x]_\ell:=\sum\limits_{{0\leq i,j\leq\kappa,\atop i+j=\ell}}x_{ij}=\pi_\ell,0\leq\ell\leq 2\kappa\},
$$
and we arrive at a semidefinite representation of $\cP$:
\begin{equation}\label{sdpP}
\cP=\{[p;q]: \exists x\in\bS^{\kappa+1}: x\succeq0\ \&\ \cQ x=P[p;q]\}.
\end{equation}
We claim that this representation is essentially strictly feasible. Indeed, let $\bar{p}=[1;1;...;1]\in\bR^{\kappa+1}$ and $\bar{q}=0$. Then
$$
\pi_{\bar{p},\bar{q}}(\tau)=(1+\tau^2)^\kappa\left[1+\tau^2/(1+\tau^2)+[\tau^2/(1+\tau^2)]^2+...+[\tau^2/(1+\tau^2)]^\kappa\right]=\sum_{i=0}^\kappa c_i\tau^{2i},\,c_i>0\,\forall i,
$$
implying that with $\bar{x}=\Diag\{c_0,c_1,...,c_\kappa\}\in\bS^{\kappa+1}$ one has
$$
P[\bar{p};\bar{q}]=\Q\bar{x}\ \&\ \bar{x}\succ0.
$$
That is, $\bar{p},\bar{q}$, and $\bar{x}$ satisfy all constraints in (\ref{sdpP}) and strictly satisfy the non-polyhedral constraint $x\succeq0$, as required by essentially strict feasibility.
\paragraph{3$^0$.} Now we are done: by its origin, $\cZ$ is a convex compact set and as such is convex and closed, implying by duality that
$$
\begin{array}{rcl}
\cZ&=&\{z\in\bR^{\kappa+1}: p^Tz-q\geq0\,\,\forall [p;q]\in\cP\}\\
&=&\Big\{z\in\bR^{\kappa+1}: \inf\limits_{p,q,x}\{[p^Tz-q]:P[p;q]-\Q x=0,x\succeq0\}{\geq 0\Big\}}\;\;\hbox{[by(\ref{sdpP})]}\\
&=&\Big\{z\in\bR^{\kappa+1}: \exists \lambda\in\bR^{2\kappa+1}: P^T\lambda=[z;-1],\cQ^*\lambda\succeq0\Big\}\\
&&\hbox{[by semidefinite duality; $\lambda\mapsto\cQ^*\lambda:\bR^{2\kappa+1}\to\bS^{\kappa+1}$}\\&&\hbox{is the conjugate of $x\mapsto\cQ x:\bS^{\kappa+1}\to\bR^{2\kappa+1}$],}
\end{array}
$$
and we arrive at the desired SDR of $\Z$.
\end{enumerate}
\subsection{Verification of calculus rules}
\begin{enumerate}
\item{\em [Restriction on a $\cK$-representable set]}
Let  (\ref{jan1}) represent $(F,\X)$. Suppose that $[t;g;x]$  can be augmented by $u,v$ to solve (\ref{jan3}). Then $(t,g,x,u)$ solve (\ref{jan1}), implying that $[t;g;x]\in\cF[F,\X]$ by definition of $\cK$-representability of $(F,\X)$, and $(x,v)$ solve (\ref{jan3}), implying that
$x\in \Y$. Taking together, these inclusions clearly imply that $[t;g;x]\in\cF[\bar{F},\Z]$, as required in item (i) of the definition of $\cK$-representation of $(\bar{F},\Z)$. Next,  when $x\in \Z$,
the triple $[t:=\langle F(x),x\rangle;g:=F(x);x]$  can be augmented by $u$ to solve (\ref{jan1}) (by item (ii) of the definition of $\cK$-representation of $(F,\X)$), and { because} $x\in\Y$,  $x$ can be augmented by $v$ to solve (\ref{jan2}). Thus, $t:=\langle F(x),x\rangle$, $g:=F(x)$, $x$ can be augmented by $(u,v)$ to solve (\ref{jan3}), as required in item (ii) of $\cK$-representability of $(\bar{F},Z)$. Thus, (\ref{jan3}) indeed represents $(\bar{F},\Z)$. A completely similar reasoning shows that if (\ref{jan1}) almost represents $(F,\X)$, then (\ref{jan3}) almost represents $(\bar{F},\Z)$.
\item{\em[Direct summation]}
When $x=[x_1;...;x_K]\in\X$ and $g=F(x)=[g_1;...;g_K]$, $g_k=F_k(x_k)$, $t_k=\langle x_k,g_k\rangle$ and $t=\sum_kt_k=\langle g,x\rangle$, there exist $u_k$ such that
(\ref{jan66}.$a$) take place (by item (ii) of the definition of representation  as applied to the representations in (\ref{jan65})), and (\ref{jan66}.$b$) takes place as well, as required in item (ii) of the definition of a representation of $(F,\X)$. On the other hand, when $x=[x_1;...;x_K]$, $g=[g_1;...;g_k]$, and $t$ can be augmented by $u_k$ and $t_k$
to solve $(\ref{jan66}.a)$ and $(\ref{jan66}.b)$, we have $[t_k;g_k;x_k]\in\cF[F_k,\X_k]$, $k\leq K$, whence for every $y=[y_1;...;y_K]\in\X$ it holds
$$
t_k-\langle g_k,y_k\rangle \geq \langle F_k(y_k),x_k-y_k\rangle,\,k\leq K.
$$
When summing up the above inequalities over $k$, we get
$$
t-\langle g,y\rangle \geq \langle F(y),x-y\rangle\,\forall y\in\X,
$$
that is, $[t;g;x]\in\cF[F,\X]$, as required in item (i) of the definition of a representation of $(F,\X)$.\\
The above reasoning, with evident modifications, shows the claim in the case of almost representations.
\item{\em[Taking conic combinations]}  Let $t,g,x$ can be augmented by $u_k,t_k,g_k,k\leq K$ to solve (\ref{jan6}). Then $[t_k;g_k;x]\in\cF[\F_k,\X]$, implying that
$$
t_k-\langle g_k,y\rangle \geq \langle F_k(y),x-y\rangle \,\,\forall y\in X;
$$
multiplying both sides by $\alpha_k$ and summing up over $k$, we get
$$
t-\langle g,y\rangle \geq \langle F(y),x-y\rangle\,\forall y\in X,
$$
that is, $[t;g;x]\in\cF[F,\X]$. On the other hand, given $x\in \X$, let us set $g_k= F_k(x)$ and $t_k=\langle F_k(x),x\rangle$. Since the $k$-th conic constraint in (\ref{jan6}) represents $(F_k,\X)$, there exist $u_k$, $k\leq K$, such that all relations (\ref{jan6}.$a$) take place. Setting $t=\sum_k\alpha_kt_k$, $g=\sum_k\alpha_kg_k$, we, on one hand, satisfy (\ref{jan6}.$b$-$c$), and, on the other hand, obtain
$$
t=\sum_k\alpha_k\langle F_k(x),x\rangle =\langle F(x),x\rangle,\,\,g=\sum_k\alpha_k g_k= F(x),
$$
the bottom line being that $[\langle F(x),x\rangle;F(x);x]$ can be augmented by $u_k$, $t_k$, $g_k$ to solve (\ref{jan6}). Thus, (\ref{jan6}) indeed is a $\cK$-representation of $(F,\X)$.\\
The above reasoning, with evident modifications, works in the case of almost representations.

\item {\em[Affine substitution of variables]}
Assume that $\tau,\gamma,\xi,g,t,u$ solve (\ref{jan8}). Then $t$, $g$, $x:=A\xi+a$, $u$ satisfy
$$
Xx+Gg+tT+Uu\leq_{\bK} a,
$$
implying that $x\in\X$ and
$$
t-\langle g,y\rangle \geq \langle F(y),x-y\rangle\,\,\forall y\in\X.
$$
Recall that $x=A\xi+a$, $x\in \X$, implies that $\xi\in\Xi$. Now, when $\eta\in\Xi$, setting $y=A\eta+a$, we have $y\in \X$ and
$$\begin{array}{rcl}
\tau-\langle\gamma,\eta\rangle&=&t-\langle g, a\rangle -\langle A^Tg,\eta\rangle =t-\langle g,y\rangle\geq\langle F(y),x-y\rangle\\
&=&
\langle F(A\eta+a),A(\xi-\eta)\rangle=\langle\Phi(\eta),\xi-\eta\rangle,\\
\end{array}
$$
and since $\eta\in\Xi$ is arbitrary, we get $[\tau;\gamma;\xi]\in  \cF[\Phi,\Xi]$, as required in item (i) of the definition of a representation of $(\Phi,\Xi)$. On the other hand, when $\xi\in \Xi$, $\gamma=\Phi(\xi)$, and $\tau=\langle \Phi(\xi),\xi\rangle$, setting
$x=A\xi+a$, we get $x\in\X$. Next, when setting $g=F(x)$, $t=\langle F(x),x\rangle$, we obtain $\gamma=A^T g$ and
\[\tau=\langle F(x),x-a\rangle=t-\langle g,a\rangle.\] Besides this, by the origin of $t,g,x$ and
item (ii) of the definition of $\cK$-representation, as applied to (\ref{jan7}), there exists $u$ such that $t,g,x,u$ satisfy (\ref{jan7}). The bottom line is that $\xi,\gamma,\tau$ can be augmented by $t,g,u$ to solve (\ref{jan8}), that is, (\ref{jan8}) meets item (ii) of definition of $\cK$-representation of $(\Phi,\Xi)$.\\
The above reasoning, with evident modifications, works in the case of almost representations.
\end{enumerate}
\section{Verifying \protect{(\ref{observe1})} and \protect{(\ref{observe2})}}
\label{justif}
{\bf A.} Let us prove that for all $u\in\U=[0,U]$, $v\in\V=[0,V]$ one has
\begin{eqnarray*}
 \psi(u,v)&=&-{u\over u+v+1}\nonumber\\
 &=&\min\limits_{f,t,s}\left\{fv+t:\, 1\geq f\geq0,\,0\leq s\leq \sqrt{uf},\,s+{1\over u+1}\leq 1,\,t-f\geq (1-s)^2-1,\,t\leq 1\right\}\nonumber\\
&=&\min\limits_{f,t,s;,\tau}\left\{fv+t:\,\begin{array}{l} 1\geq f\geq0,\,s\geq0,\,t\leq 1,\,\tau\geq0,\,s+\tau\leq 1\\
\underbrace{[2s;u-f;u+f]\in\bL^3}_{\Leftrightarrow s^2\leq uf\hbox{\tiny\ when $u+f\geq0$}}\\
\underbrace{[2(1-s);t-f;t-f+2]\in\bL^3}_{\Leftrightarrow (1-s)^2\leq t-f+1\hbox{\tiny\ when $t-f+2\geq0$}}\\
\underbrace{[2;u-\tau+1;u+\tau+1]\in\bL^3}_{\Leftrightarrow (u+1)\tau\geq1\hbox{\tiny\ when $u+\tau+1\geq0$}}
\end{array}\right\}.
\end{eqnarray*}
Indeed, for $u\in\cU$ and $v\in\cV$ we have
$$
\begin{array}{l}
\min\limits_{f,t,s}\left\{fv+t:\, 1\geq f\geq0,0\leq s\leq \sqrt{uf},s+{1\over u+1}\leq 1,t-f\geq (1-s)^2-1,t\leq 1\right\}\\
=\min\limits_{f}\left\{fv+f-2\bar{s}(f)+\bar{s}(f)^2:\,1\geq f\geq0,\,\bar{s}(f)=\min\left[{u\over u+1},\sqrt{uf}\right]\right\}\\
=\min\left[\min\limits_{0\leq f\leq u/(u+1)^2}\left\{fv+f-2\sqrt{uf}+uf\right\},\min\limits_{1\geq f\geq u/(u+1)^2}
\left\{fv+f-{2u\over u+1}+{u^2\over(u+1)^2}\right\}\right]\\
=\min\left[-{u\over u+v+1},{u(v+1)-2u-u^2\over (u+1)^2}\right]=u\min\left[-{1\over v+(u+1)},{v-(u+1)\over (u+1)^2}\right]={u\over v+u+1}\min\left[-1,{v^2-(u+1)^2\over(u+1)^2}\right]\\
=\psi(u,v),\\
\end{array}
$$
as claimed.
\par\noindent
{B.} Now let us show that for all $u\in\U=[0,U]$, $v\in\V=[0,V]$
\begin{eqnarray*}
 -\psi(u,v)&:=&{u\over u+v+1}\\
 &=&\min\limits_{f,t,s}\left\{fu+t: 1\geq f\geq0,\,0\leq s\leq \sqrt{(v+1)f},\,s\leq1,\,1\geq t\geq (1-s)^2\right\}\\
&=&\min\limits_{f,t,s}\left\{fu+t:\begin{array}{l} 1\geq f\geq0,\,0\leq s\leq 1,\,t\leq 1\\
\underbrace{[2s;v+1-f;v+1+f]\in\bL^3}_{\Leftrightarrow s^2\leq (v+1)f\hbox{\tiny\ when $v+1+f\geq0$}}\\
\underbrace{[2(1-s);t-1;t+1]\in\bL^3}_{\Leftrightarrow (1-s)^2\leq t}\\
\end{array}\right\}
\end{eqnarray*}
Indeed, for $u\in\cU$ and $v\in\cQ$ we have
$$
\begin{array}{l}
\min\limits_{f,t,s}\left\{fu+t: 1\geq f\geq0,0\leq s\leq \sqrt{(v+1)f},s\leq 1,1\geq t\geq (1-s)^2\right\}\\
=
\min\limits_{f,s}\left\{fu+(1-s)^2: 1\geq f\geq0,0\leq s\leq \sqrt{(v+1)f},s\leq 1\right\}\\
=\min\limits_{f}\left\{fu+(1-\bar{s}(f))^2: 1\geq f\geq0,\bar{s}(f)=\min[\sqrt{(v+1)f},1]\right\}\\
=\min\left[\min\limits_{0\leq f\leq 1/(v+1)}\{fu+1-2\sqrt{(v+1)f}+(v+1)f\},\,\min\limits_{1\geq f\geq 1/(v+1)}\{fu\}\right]\\
=\min\left[{u\over u+v+1}, {u\over v+1}\right]={u\over u+v+1}=-\psi(u,v),\\
\end{array}
$$
as claimed.
\qed
 \end{document}